\documentclass[10pt]{amsart}
\usepackage{amsmath,amstext,amsfonts,%
enumitem,marvosym,%
extpfeil,graphicx,mathrsfs,%
mathtools,setspace,tikz,xparse}
\usepackage[alphabetic]{amsrefs}
\usepackage[unicode]{hyperref}
\usepackage[utf8]{inputenc}
\usepackage[T1]{fontenc}
\usepackage[normalem]{ulem}
\usepackage[cmtip,all]{xy}
\usepackage[light]{sourcesanspro}
\usepackage{lmodern}
\hypersetup{colorlinks=false,linkbordercolor={0.93 0.71 0.13},citebordercolor={1 0 .5},urlbordercolor={0.27 0.39 0.68}}%
\usetikzlibrary{decorations.text,shapes.geometric,arrows,positioning,fit,calc,matrix}

\usepackage{filecontents}

\definecolor{MyColor1}{HTML}{fff3b0}
\definecolor{MyColor2}{HTML}{e09f3e}
\definecolor{MyColor3}{HTML}{335c67}
\definecolor{MyColor4}{HTML}{393e41}

\DeclareMathOperator{\GL}{GL}

\DeclareMathOperator{\Img}{im}

\DeclareMathOperator{\ind}{ind}

\newcommand{\m}[1]{\(\smash{#1}\)}

\newcommand{\bb}[1]{\mathbb{#1}}
\newcommand{\QQ}{\bb{Q}}
\newcommand{\FF}{\bb{F}}
\newcommand{\ZZ}{\bb{Z}}
\newcommand{\CC}{\bb{C}}
\newcommand{\FFp}{\overline{\bb{F}}_p}
\newcommand{\QQp}{\smash{\overline{\bb{Q}}_p}}

\newcommand{\ZZp}{\smash{\overline{\bb{Z}}_p}}

\newcommand{\defeq}{\vcentcolon=}
\newcommand{\SCR}{\mathscr}
\newcommand{\IE}{i.e.\ }

\newcommand{\HECKET}{\SCR T}

\newcommand{\MU}{\mu}

\newcommand{\SIGMA}{\Sigma}

\newcommand{\LAMBDA}{\lambda}

\newcommand{\LESS}[1]{\smash{\uline{#1}}}
\newcommand{\MORE}[1]{\smash{\overline{#1}}}
\newcommand{\BIGO}[1]{\smash{\mathsf{O}\!\left({#1}\right)}}

\newcommand{\SIZEEE}{4.25pt}

\newcommand{\IRREDUCIBLE}{\mathrm{Irr}}
\newcommand{\RREDUCIBLE}{\mathrm{Red}}
\newcommand{\BIRREDUCIBLE}{\mathrm{BIrr}}
\newcommand{\BRREDUCIBLE}{\mathrm{BRed}}
\newcommand{\TEMPORARYA}{\xi}

\newcommand{\TEMPORARYD}{\vartheta}%

\newcommand{\TEMPORARYI}{\nu}
\newcommand{\TEMPORARYM}{\mathrm{M}}

\newcommand{\PADICVALUATION}{v_p}

\newcommand{\TEMPORARYF}{L}
\newcommand{\TEMPORARYDISC}{D}

\newcommand{\sfirst}{s_{1}}
\newcommand{\ssecond}{s_{2}}
\newcommand{\TILDE}{\widetilde}

\newcommand{\GEQSLANT}{\geqslant}
\newcommand{\LEQSLANT}{\leqslant}

\newcommand{\maxideal}{\mathfrak m}
\newcommand{\repnn}{\pi}
\newcommand{\NEWLINEA}{\\[2pt]}

\newcommand{\SUBMODULE}{\mathrm{sub}}
\newcommand{\QUOTIENTI}{\mathrm{quot}}

\newcommand{\FILT}{\widehat}
\newcommand{\FORMX}{X}
\newcommand{\FORMY}{Y}
\newcommand{\DELTAA}{\Delta}
\newcommand{\TEMPORARYJ}{\eta}

\newcommand{\sqbrackets}[4]{{#1}\bullet_{#2,#3}{#4}}

\newtheorem{theorem}{Theorem}
\newtheorem{lemma}[theorem]{Lemma}
\newtheorem{proposition}[theorem]{Proposition}
\newtheorem{corollary}[theorem]{Corollary}
\newtheorem{conjecture}{Conjecture}

\newenvironment{Smatrix}%
    {\tiny (\begin{smallmatrix}}%
    {\end{smallmatrix})}
\newenvironment{Lmatrix}%
    {\footnotesize (\begin{smallmatrix}}%
    {\end{smallmatrix})}
\newcommand{\BLACKSQUARE}{\hspace{\stretch{1}}\m{\mathbin{\vcenter{%
    \hbox{\rule{.9ex}{.9ex}}}}}}
\renewenvironment{proof}[1]%
    {\emph{{#1}}}%
    {\BLACKSQUARE}

\makeatletter
\newcommand{\xhookdoubleheadrightarrow}[2][]{%
  \lhook\joinrel
  \ext@arrow 0359\rightarrowfill@ {#1}{#2}%
  \mathrel{\mspace{-15mu}}\rightarrow}
\renewcommand{\xtwoheadrightarrow}[2][]{%
  \ext@arrow 0359\rightarrowfill@ {#1}{#2}%
  \mathrel{\mspace{-15mu}}\rightarrow}
\def\@settitle{\begin{center}%
  \baselineskip14\p@\relax\normalfont\LARGE
    \uppercasenonmath\@title
  \@title\end{center}%
}
\def\do#1{\@ifdefinable#1{\newdimen#1}}
\do\@TempWidthOne
\do\@TempWidthTwo
\newcommand*\WidthInMathUnitsOf[2]{%
    \settowidth\@TempWidthOne{$\m@th #1#2$}%
    \settowidth\@TempWidthTwo{$\m@th #1\mspace{1mu}$}%
    \strip@pt \dimexpr \@TempWidthOne*\p@/\@TempWidthTwo \relax
}
\newcommand*\SetToWidthInMathUnits[3]{%
    \settowidth\@TempWidthOne{$\m@th #2#3$}%
    \settowidth\@TempWidthTwo{$\m@th #2\mspace{1mu}$}%
    \edef #1{%
        \strip@pt \dimexpr \@TempWidthOne*\p@/\@TempWidthTwo \relax
    }%
}
\newcommand*\ExtractHexDigitOfFamilyOfSymbol[1]{%
    \expandafter\@ExtractHexDigitOfFamily \meaning #1%
}
\@ifdefinable\@ExtractHexDigitOfFamily{
    \def\@ExtractHexDigitOfFamily#1"#2#3#4#5{#3}
}
\makeatother

\newcommand*{\ReportValueOfOneEm}[1]{%
    &\text{in #1 size,}%
        &\quad\texttt{1em}%
            &=\texttt{\the \fontdimen 6 \csname #1font\endcsname 2}%
        &&\qquad\text{(from \texttt{\fontname \csname #1font\endcsname 2})}%
}
\newcommand*{\ReportSizeOfFontOfSymbol}[2]{%
    &\texttt{\pdffontsize\csname #2font\endcsname
                "\ExtractHexDigitOfFamilyOfSymbol{#1}}%
        &&\quad\text{in #2 size}%
        &\qquad(\texttt{\string\fontname}%
            &=\texttt{\fontname\csname #2font\endcsname
                "\ExtractHexDigitOfFamilyOfSymbol{#1}})%
}

\title%
[Reductions of crystalline representations, II]%
{On the reductions of certain
two-dimensional crystalline representations, II}
\author{Bodan Arsovski}
\date{February 2021}

\begin{document}
\maketitle

\renewcommand{\HECKET}{T}
\newcommand{\RRRRRR}{\SCR R^{s,\TEMPORARYI}}
\newcommand{\VERIFYIF}[1]{[{#1}]}
\newcommand{\DARKCIRC}{\bullet}
\newcommand{\TEMPORARYAAAA}{w}
\newcommand{\UUU}{U}
\newcommand{\subm}{sub}

\newcommand{\MATRIXNO}{\SCR M_1}
\newcommand{\MATRIXYES}{\SCR M_2}
\renewcommand{\TEMPORARYF}{m}
\renewcommand{\TEMPORARYDISC}{\SCR N}
\newcommand{\LARGESTSLOPE}{100}
\newcommand{\LARGESTPRIME}{139522597633}
\newcommand{\ddd}{\mathrm{d}}
\newcommand{\polyr}{z}
\newcommand{\polyH}{\Phi}
\newcommand{\polys}{\psi}
\newcommand{\polyconst}{C}
\newcommand{\veeC}{\check{C}}
\newcommand{\vveeC}{\hat{C}}
\newcommand{\annd}{\text{ \& }}
\newcommand{\Cnegativeone}{\frac{(-1)^{\beta}(s-\alpha)(\alpha-\beta+1)}{\beta^2(2\alpha-s+1) \binom{\alpha}{\beta}} \binom{\alpha}{s-\alpha}}
\newcommand{\rowR}{\mathbf{r}}
\newcommand{\colA}{\mathbf{c}}
\newcommand{\abcNL}{\\[3pt]}
\newcommand{\Cjnewline}{\\[5pt]}
\newcommand{\tempheightaer}{\rule{0pt}{12pt}}
\newcounter{equ}
\newcommand{\jj}{t}
\newcommand{\llll}{\gamma}
\newcommand{\zs}{u}
\newcommand{\za}{v}
\newcommand{\zb}{w}
\newcommand{\zi}{t}
\newcommand{\varu}{e}
\newcommand{\boxcol}[1]{\tikz{\path[draw=black,fill={#1}] (0,0) rectangle (0.2cm,0.2cm);}}

\begin{abstract}
A conjecture of Breuil, Buzzard, and Emerton says that the slopes %
 of certain  %
 reducible \m{p}-adic Galois representations must be integers. In~\cite{Ars1} we showed this conjecture for %
  representations that lie over certain ``non-subtle'' components of weight space. %
   This article is a continuation of~\cite{Ars1} in which we completely  classify the aforementioned representations over the ``non-subtle'' components of weight space, both for integer and non-integer slopes.
\end{abstract}

\section{Introduction and results}

\subsection{Background} Let \m{p} be an odd prime number %
 and \m{k\GEQSLANT2} be an integer, and let \m{a} be an element of  \m{\ZZp} such that \m{\PADICVALUATION(a)>0}. Let us denote \m{\TEMPORARYI = \lfloor \PADICVALUATION(a) \rfloor + 1 \in \ZZ_{>0}}. With this data one can associate a certain two-dimensional crystalline \m{p}-adic representation \m{V_{k,a}} with Hodge--Tate weights \m{(0,k-1)}. We give the definition of this representation in section~2 of~\cite{Ars1}, and we define \m{\overline{V}_{k,a}} as the semi-simplification of the reduction modulo the maximal ideal \m{\maxideal} of \m{\ZZp} of a Galois stable \m{\ZZp}-lattice in \m{{V_{k,a}}} (with the resulting representation being independent of the choice of lattice). The question of computing \m{\overline{V}_{k,a}} has been studied extensively, and we refer to the introduction of~\cite{Ars1} for a brief exposition of it. Partial results have been obtained by %
  Fontaine, Edixhoven, Breuil, Berger, Li, Zhu, Buzzard, Gee, Bhattacharya, Ganguli, Ghate, et al (see~\cite{Ber10}, \cite{Bre03a}, \cite{Bre03b}, \cite{Edi92}, \cite{BLZ04}, \cite{BG15}, \cite{BG09}, \cite{BG13}, \cite{GG15}). 
A conjecture of Breuil, Buzzard, and Emerton says the following.
\begin{conjecture}\label{MAINCONJ}
{If \m{k} is even and \m{\PADICVALUATION(a) \not\in \ZZ} then \m{\overline{V}_{k,a}} is irreducible.}
\end{conjecture}

The main result of~\cite{Ars1} is that this conjecture is true over
certain ``non-subtle'' components of weight space.
We say that a weight \m{k} belongs to a  ``non-subtle'' component of weight space if and only if
\[\textstyle \textstyle k \not\equiv 3,4,\ldots,2\TEMPORARYI,2\TEMPORARYI+1 \bmod p-1.\]
Thus there are \m{\max\{\frac{p-1}{2} - \TEMPORARYI+1,0\}} many ``non-subtle'' components of weight space. 
This article is a continuation of~\cite{Ars1} in which we completely classify these representations over the ``non-subtle'' components of weight space, both for integer and non-integer slopes. 
Before stating the main results, let us first introduce some terminology.

We fix embeddings \m{\overline{\QQ}\hookrightarrow \QQp\hookrightarrow \CC_p} and  
  to look at the space of continuous homomorphisms \m{\SCR W = \hom_{cts}(\ZZ_p^\times,\CC_p^\times)}. We have \m{\ZZ_p^\times \cong (\ZZ/(p-1)\ZZ) \times \ZZ_p}, and \m{\hom_{cts}(\ZZ_p,\CC_p^\times)} is isomorphic to the open disk in \m{\CC_p} with center \m{0} and radius \m{1}, via the identification \m{\chi \leftrightarrow \chi(1)-1}.
Thus \m{\SCR W} is the disjoint union of \m{p-1} open disks of radii 1. %
 We
 identify the weight \m{k\GEQSLANT 2} with the continuous homomorphism \m{x \mapsto x^{k-2}}, so that \m{k\GEQSLANT 2} is a point on the disk indexed by \m{k-2 \bmod {p-1}}. We say a point of \m{\SCR W} is ``integral'' if it is associated with a weight in this fashion.
 
 The most general result about \m{\overline{V}_{k,a}} to date is the main theorem in~\cite{BLZ04} which says that \[\textstyle\overline{V}_{k,a} \cong \overline{V}_{k,0}\cong\left\{\begin{array}{rl} %
\ind(\omega_2^{k-1}) & \text{if } p+1\nmid k-1, \\[2pt] (\mu_{\sqrt{-1}} \oplus \mu_{-\sqrt{-1}})\otimes \omega^{\frac{k-1}{p+1}} & \text{if } p+1\mid k-1,\end{array}\right.\] whenever \m{\PADICVALUATION(a) > \lfloor\frac{k-2}{p-1}\rfloor}. This theorem tells us what \m{\overline{V}_{k,a}} is at a discrete set of points in the aforementioned \m{p-1} disks. It is expected that the bound \m{\lfloor\frac{k-2}{p-1}\rfloor} is not optimal: there is a prediction in subsection~2.1 of~\cite{BG16} that the bound can be replaced with \m{\frac{k-1}{p+1}}.

The theorems we prove indicate that these points play a fundamental role: it seems that the \m{p-1} open disks can be split into regions by concentric circles centered at these points in a way that \m{\overline{V}_{k,a}} depends on the region \m{k} belongs to.
For example, the following diagram  illustrates how the disk containing the weight \m{10} is split into regions when \m{p\GEQSLANT 11} and \m{3 < \PADICVALUATION(a) < 4}. There are four centers of the bundles of concentric circles, and they are exactly the points \m{k} belonging to the disk (\IE such that \m{k\equiv 10 \bmod p-1}) and such that \m{\PADICVALUATION(a) > \lfloor\frac{k-2}{p-1}\rfloor} (\IE \m{\lfloor\frac{k-2}{p-1}\rfloor < 4}). In fact, they are the same as the points satisfying \m{\frac{k-1}{p+1} < 4}. %
By~\cite{BLZ04}, at these four points we have \m{\overline{V}_{k,a} \cong \ind(\omega_2^{k-1})}.
The diagram illustrates closed disks around the points, of radii \m{p^{-1},p^{-2}}, and \m{p^{-3}}.  As we show in theorem~\ref{MAINTHM0}, in this situation \m{\overline{V}_{k,a}} is always one of these four representations,  depending on the color of the region that \m{k} belongs to, as illustrated.

\null

 \vspace{\stretch{1}}

\begin{center}
\begin{tikzpicture}[
        scale=0.6,
        black,
        ultra thick,
        planet/.style = {draw,fill,circle,inner sep=#1},
        circle label/.style = {
            postaction={
                decoration={
                    text along path,
                    text = {#1},
                    text align=center,
                    text color=black,
                    reverse path,
                },
            decorate,
        }
        }
    ]
    \foreach \i in {6} {
        \draw[color=white,fill=MyColor1] (0,0) circle (\i);
    }
    \draw[color=MyColor2,fill=MyColor2] (-3.5,0) circle (2.1);
    \draw[color=MyColor3,fill=MyColor3] (-3.5,0) circle (1.7);
    \draw[color=MyColor4,fill=MyColor4] (-3.5,0) circle (1.3);
    \draw[color=MyColor2,fill=MyColor2] (3.5,0) circle (2.1);
    \draw[color=MyColor2,fill=MyColor2] (0,3.5) circle (2.1);
    \draw[color=MyColor3,fill=MyColor3] (0,3.5) circle (1.7);
    \node[white,planet=0.5pt,label={\footnotesize\bf {\color{white}${10}$}%
    }] at (-3.5,0) {};
    \node[white,planet=0.5pt,label={\footnotesize\bf {\color{white}${p+9}$}%
    }] at (0,3.5) {};
    \node[black,planet=0.5pt,label={\footnotesize\bf ${2p+8}$%
    }] at (3.5,0) {};
    \node[black,planet=0.5pt,label={\footnotesize\bf ${3p+7}$%
    }] at (0,-3.5) {};
\end{tikzpicture}
\end{center}

 Theorems~\ref{MAINTHM0} and~\ref{MAINTHM00} prove a general version of this, in the case when the corresponding disk is what we label ``non-subtle'': this label depends on the slope \m{\PADICVALUATION(a)} and roughly means that the bound \m{\lfloor\frac{k-2}{p-1}\rfloor} is optimal for any \m{k} belonging to that disk, \IE that there are no additional points \m{k} on that disk satisfying the improved bound \m{\PADICVALUATION(a) > \frac{k-1}{p+1}} but not \m{\PADICVALUATION(a) > \lfloor\frac{k-2}{p-1}\rfloor} at which the associated representation is distinct from the ``typical'' one, which in the context of the diagram on the previous page means the representation corresponding to the yellow region.
 
 Let us write 
 \m{\MORE{h}} for the number in \m{\{1,\ldots,p-1\}} which is congruent to \m{h} mod \m{p-1}.  Let \m{\TEMPORARYI = \lfloor \PADICVALUATION(a) \rfloor + 1 \in \ZZ_{>0}}, and let \m{s} be the number in \m{\{1,\ldots,p-1\}} which is congruent to \m{k-2} mod \m{p-1}. Let us say that \m{k} is ``subtle'' if  \m{s \in\{1,\ldots,2\TEMPORARYI-1\}}, and \m{k} is ``non-subtle'' if  \m{s \not\in\{1,\ldots,2\TEMPORARYI-1\}}.
 An open disk of \m{\SCR W} consists either entirely of ``subtle'' points or entirely of ``non-subtle'' points, so we can also refer to the \m{p-1} open disks of \m{\SCR W} as either ``subtle'' or ``non-subtle''. In particular,  the \m{\min\{p-1,2\TEMPORARYI-1\}} disks containing \m{3,\ldots,2\TEMPORARYI+1} are ``subtle'', and all other disks are ``non-subtle''. Note that whether a weight is ``subtle'' or not depends on the value of \m{\TEMPORARYI}.

  Let \m{\SCR D_s} denote the open disk of radius 1 around \m{s+2\in \SCR W}, and  let us consider the set
\[\textstyle B_{s,\TEMPORARYI} = \textstyle \{s + \beta(p-1)+2 \,|\, \beta \in \{0,\ldots,\TEMPORARYI-2\}\}.\]
     In particular, if \m{\TEMPORARYI=1} then \m{B_{s,\TEMPORARYI}=\varnothing}, and in general \m{B_{s,\TEMPORARYI}} is a set of \m{\TEMPORARYI-1} points in \m{\SCR D_s}. Let \[b_0 < \cdots < b_{\TEMPORARYI-2}\] be the elements of \m{B_{s,\TEMPORARYI}} in increasing order. Therefore if \m{i \in \{0,\ldots,\TEMPORARYI-2\}} then \m{b_i = s+i(p-1)+2} and, by the main result of~\cite{BLZ04}, \[\textstyle \overline V_{b_i,a}\cong \ind(\omega_2^{{b_i-1}}).\] Let us also define \m{b_{\TEMPORARYI-1} = s+(\TEMPORARYI-1)(p-1)+2}. For \m{i \in \{0,\ldots,\TEMPORARYI-2\}} and \m{j \in \ZZ_{>0}}, let \[\textstyle\RRRRRR_{i,j} = \{t \in \SCR D_s \,|\, j\LEQSLANT \PADICVALUATION(t-b_i) < j+1\}\] be the half-open annulus which is the complement of the closed disk of radius \m{p^{-j-1}} around \m{b_i} in the closed disk of radius \m{p^{-j}} around \m{b_i}. The integral points in \m{\RRRRRR_{i,j}} are the points on the circle of radius \m{p^{-j}} around \m{b_i = s+i(p-1)+2}. Finally, let \[\textstyle \RRRRRR_0 = \SCR D_s \backslash \cup_{i\in \{0,\ldots,\TEMPORARYI-2\}, \, j > 0} \RRRRRR_{i,j},\]  so that \m{\SCR D_s} is partitioned into the disjoint sets \begin{align}\label{disjset}\tag{\m{\RRRRRR}} \textstyle \{\RRRRRR_0\} \cup \{\RRRRRR_{i,j}\,|\, i \in \{0,\ldots,\TEMPORARYI-2\},\, j\in\ZZ_{>0}\}.\end{align}
Note that the definition of this partition depends on both \m{s} and \m{\TEMPORARYI}. For \m{l \in \ZZ} and \m{\lambda\in \smash{\overline{\FF}}_p^\times} let us define
       \[\textstyle  \IRREDUCIBLE(l) = \ind(\omega_2^{l-1}) \text{ and }\RREDUCIBLE_{s,\TEMPORARYI}(l,\lambda) = \mu_{\lambda} \omega^{s+l-\TEMPORARYI+2} \oplus \mu_{\lambda^{-1}} \omega^{\TEMPORARYI-l-1}.\]
The first result %
 is a complete classification of %
 \m{\overline V_{k,a}} over the ``non-subtle'' components of weight space for \m{\PADICVALUATION(a) \not\in \ZZ}.
     
     \begin{theorem}\label{MAINTHM0}
Recall that \m{k\GEQSLANT 2} is an integer and that \m{s} is defined as the integer in \m{\{1,\ldots,p-1\}} which is congruent to \m{k-2 \bmod p-1}. Suppose that \m{k} is ``non-subtle'', \IE 
 \[\textstyle k \not\equiv 3,4,\ldots,2\TEMPORARYI,2\TEMPORARYI+1 \bmod p-1.\]
 Suppose also that the open disk \m{\SCR D_s} of radius 1 around \m{s+2\in \SCR W} is partitioned into disjoint sets as in~\eqref{disjset}. %
 If \m{\PADICVALUATION(a)\not\in\ZZ} then
\[\textstyle \overline{V}_{k,a} \cong \left\{\begin{array}{rl} %
\IRREDUCIBLE(b_{\TEMPORARYI-1})
& \text{if } k \in \RRRRRR_0, \\[2pt] \IRREDUCIBLE(b_{\max\{i,\TEMPORARYI-j-1\}})  & \text{if } k \in \RRRRRR_{i,j}.\end{array}\right.\]
\end{theorem}

This result is known for \m{\TEMPORARYI=1} by the work of Buzzard and Gee in~\cite{BG09} and for \m{\TEMPORARYI=2} by the work of Bhattacharya and Ghate in~\cite{BG15}. %

 We also prove a similar theorem for \m{\PADICVALUATION(a)\in\ZZ}. The precise statement of this is the following theorem, which is a complete classification of %
 \m{\overline V_{k,a}} over the ``non-subtle'' components of weight space for \m{\PADICVALUATION(a) \in \ZZ}.
\begin{theorem}\label{MAINTHM00}
Recall that \m{k\GEQSLANT 2} is an integer and that \m{s} is defined as the integer in \m{\{1,\ldots,p-1\}} which is congruent to \m{k-2 \bmod p-1}. Suppose that \m{k} is ``non-subtle'', \IE 
 \[\textstyle k \not\equiv 3,4,\ldots,2\TEMPORARYI,2\TEMPORARYI+1 \bmod p-1.\]
 Suppose also that the open disk \m{\SCR D_s} of radius 1 around \m{s+2\in \SCR W} is partitioned into disjoint sets as in~\eqref{disjset}.  %
 If \m{\PADICVALUATION(a)=\TEMPORARYI-1\in\ZZ_{>0}} then
\[\textstyle \overline{V}_{k,a} \cong \left\{\begin{array}{rl} %
\RREDUCIBLE_{s,\TEMPORARYI}(0,\lambda_{k,\TEMPORARYI}) & \text{if } k \in \RRRRRR_0, \\[2pt] \RREDUCIBLE_{s,\TEMPORARYI}(j,\lambda_{k,\TEMPORARYI,i,j}) & \text{if } k \in \RRRRRR_{i,j} \text{ and } i+j < \TEMPORARYI-1,\\[2pt]
\IRREDUCIBLE(b_i) & \text{if } k \in \RRRRRR_{i,j} \text{ and } i+j \GEQSLANT \TEMPORARYI-1,\end{array}\right.\]
where 
\begin{align*}
    \textstyle \lambda_{k,\TEMPORARYI} & \textstyle = \frac{\binom{s-\TEMPORARYI+2}{\TEMPORARYI-1}a}{\binom{s-k+2}{\TEMPORARYI-1}p^{\TEMPORARYI-1}} \in \smash{\overline{\FF}}_p^\times,
    \\ \textstyle \lambda_{k,\TEMPORARYI,i,j} & \textstyle = \frac{(-1)^{\TEMPORARYI+i+j+1} (\TEMPORARYI-j-1) \binom{\TEMPORARYI-j-2}{i}\binom{s-\TEMPORARYI+j+2}{\TEMPORARYI-j-1}a}{(k-s-i(p-1)-2)p^{\TEMPORARYI-j-1}} \in \smash{\overline{\FF}}_p^\times.
\end{align*}
\end{theorem} 
 Note that \m{\lambda_{k,\TEMPORARYI,i,j}} is indeed a unit (and therefore we can think of it as an element of \m{\smash{\overline{\FF}}_p^\times}) since the integral points in  \m{\RRRRRR_{i,j}} consist precisely of those points on the circle of radius \m{p^{-j}} around \m{b_i = s+i(p-1)+2} and so 
\[\textstyle \PADICVALUATION((k-s-i(p-1)-2)p^{\TEMPORARYI-j-1}) = \TEMPORARYI-1.\]
This result is known for \m{\TEMPORARYI=2} by the work of Bhattacharya, Ghate, and Rozensztajn in~\cite{BGR18}. 

The paper~\cite{Roz18} gives an algorithm which takes as input a prime \m{p}, a weight \m{k}, an eigenvalue \m{a}, and a parameter called ``radius'' which determines the precision of the computations, and if the radius is large enough it computes \m{\overline V_{k,a}}.\footnote{The algorithm actually computes the  \m{\GL_2(\QQ_p)}-representation associated with \m{\overline V_{k,a}} via the bijective correspondence given in Theorem~2 in~\cite{Ars1}, and one can use Theorem~2 in~\cite{Ars1} to then compute \m{\overline V_{k,a}}.} An implementation of the algorithm in \texttt{SageMath} is available at
\begin{center}\footnotesize
    \href{http://perso.ens-lyon.fr/sandra.rozensztajn/software/index.html}{\texttt{http://perso.ens-lyon.fr/sandra.rozensztajn/software/index.html}}
\end{center}
As theorems~\ref{MAINTHM0} and~\ref{MAINTHM00} give complete classifications of %
 \m{\overline V_{k,a}}, one can use this algorithm to verify them for any given triple \m{(p,k,a)}. The complexity of the algorithm depends on the size of the extension field generated by \m{a}, so in practice it is much faster to verify theorem~\ref{MAINTHM00}. Additionally, the statement of theorem~\ref{MAINTHM00} is more complicated, especially the formulas for \m{\lambda_{k,\TEMPORARYI},\lambda_{k,\TEMPORARYI,i,j}}, so it is better suited for this type of computer verification. We have verified theorem~\ref{MAINTHM00} for the triples in the following table.
 
 \
 
 \
 
 \begin{center}
     \begin{tabular}{cccc|l}
        \m{p} & \m{k} & \m{a} & ``radius'' & \ \m{\overline V_{k,a}} \\[4pt]\hline
        7 & 8 & 49 & 3 & \ \m{\ind(\omega_2^7)}\rule{0pt}{16pt} \\[4pt] %
        7 & 14 & 49 & 3 & \ \m{\ind(\omega_2^{13})} \\[4pt] %
        7 & 20 & 49 & 4 & \ \m{\mu_3 \omega^5 \oplus \mu_5 \omega^2} \\[4pt] %
        7 & 26 & 49 & 3 & \ \m{\mu_1 \omega^5 \oplus \mu_1 \omega^2} \\[4pt]%
        7 & 32 & 49 & 3 & \ \m{\mu_4 \omega^5 \oplus \mu_2 \omega^2} \\[4pt]%
        7 & 38 & 49 & 3 & \ \m{\mu_1 \omega^5 \oplus \mu_1 \omega^2} \\[4pt]%
        7 & 44 & 49 & 3 & \ \m{\mu_3 \omega^5 \oplus \mu_5 \omega^2} \\[4pt]%
        7 & 50 & 49 & 3 & \ \m{\mu_6 \omega \oplus \mu_6}\\[4pt]%
        7 & 56 & 49 & 3 & \ \m{\ind(\omega_2^{13})}\\[4pt]%
        11 & 38 & 121 & 3 & \ \m{\mu_7 \omega^5 \oplus \mu_8 \omega^2} \\[4pt]%
        11 & 39 & 121 & 3 & \ \m{\mu_5 \omega^6 \oplus \mu_9 \omega^2} \\[4pt]%
        11 & 40 & 121 & 3 & \ \m{\mu_7 \omega^7 \oplus \mu_8 \omega^2} \\[4pt]%
        11 & 41 & 121 & 3 & \ \m{\mu_2 \omega^8 \oplus \mu_6 \omega^2} \\[4pt]%
        11 & 42 & 121 & 3 & \ \m{\mu_1 \omega^9 \oplus \mu_1 \omega^2} \\[4pt]%
     \end{tabular}
 \end{center}

    \section{Computing \texorpdfstring{\m{\overline V_{k,a}}}{V} by computing \texorpdfstring{\m{\overline\Theta_{k,a}}}{\textTheta}}
    \label{secComp}
    
    From now on we assume the notation from sections 2, 3, 4, and 6 of~\cite{Ars1}. Moreover, we assume that \m{k>p^{100}} as in section~5 of~\cite{Ars1}.
         For \m{l \in \ZZ}
         let us define
       \[\textstyle  \BIRREDUCIBLE(l) = \left(\ind_{K Z}^G \sigma_{l_1}/{T}\right) \otimes \omega^{l_2},\]
       where \m{l_1} and \m{l_2} are the unique integers such that
        \m{l = l_1 + (p+1)l_2 + 2} and \m{l_1 \in \{0,\ldots,p-1\}}.
        For \m{l \in \ZZ}
         and \m{\lambda\in \smash{\overline{\FF}}_p^\times}  let us define \begin{align*} & \textstyle \BRREDUCIBLE_{s,\TEMPORARYI}(l,\lambda) \\ & \textstyle \null\qquad=
        \repnn(\LESS{s+2l-2\TEMPORARYI+2},\lambda,\omega^{\TEMPORARYI-l-1}) \oplus\repnn(\LESS{2\TEMPORARYI-s-2l-4},\lambda^{-1},\omega^{s+l-\TEMPORARYI+2}) .
\end{align*} Theorem~2 in~\cite{Ars1} implies that our main theorems can be rewritten in the following equivalent forms. Recall that we assume \m{p>2} throughout.
    
     \begin{theorem}\label{MAINTHM0a}
Recall that \m{k\GEQSLANT 2} is an integer and that \m{s} is defined as the integer in \m{\{1,\ldots,p-1\}} which is congruent to \m{k-2 \bmod p-1}. Suppose that \m{k} is ``non-subtle'', \IE 
 \[\textstyle k \not\equiv 3,4,\ldots,2\TEMPORARYI,2\TEMPORARYI+1 \bmod p-1.\]
 Suppose also that the open disk \m{\SCR D_s} of radius 1 around \m{s+2\in \SCR W} is partitioned into disjoint sets as in~\eqref{disjset}. 
 If \m{\PADICVALUATION(a)\not\in\ZZ} then
\[\textstyle \overline{\Theta}_{k,a} \cong \left\{\begin{array}{rl} 
\BIRREDUCIBLE(b_{\TEMPORARYI-1})
& \text{if } k \in \RRRRRR_0, \\[2pt] \BIRREDUCIBLE(b_{\max\{i,\TEMPORARYI-j-1\}})  & \text{if } k \in \RRRRRR_{i,j}.\end{array}\right.\]
\end{theorem}

\begin{theorem}\label{MAINTHM00a}
Recall that \m{k\GEQSLANT 2} is an integer and that \m{s} is defined as the integer in \m{\{1,\ldots,p-1\}} which is congruent to \m{k-2 \bmod p-1}. Suppose that \m{k} is ``non-subtle'', \IE 
 \[\textstyle k \not\equiv 3,4,\ldots,2\TEMPORARYI,2\TEMPORARYI+1 \bmod p-1.\]
 Suppose also that the open disk \m{\SCR D_s} of radius 1 around \m{s+2\in \SCR W} is partitioned into disjoint sets as in~\eqref{disjset}. 
 If \m{\PADICVALUATION(a)=\TEMPORARYI-1\in\ZZ_{>0}} then
\[\textstyle \overline{\Theta}_{k,a}^{\rm ss} \cong \left\{\begin{array}{rl} 
\BRREDUCIBLE_{s,\TEMPORARYI}(0,\lambda_{k,\TEMPORARYI})^{\rm ss} & \text{if } k \in \RRRRRR_0, \\[2pt] \BRREDUCIBLE_{s,\TEMPORARYI}(j,\lambda_{k,\TEMPORARYI,i,j})^{\rm ss} & \text{if } k \in \RRRRRR_{i,j} \text{ and } i+j < \TEMPORARYI-1,\\[2pt]
\BIRREDUCIBLE(b_i) & \text{if } k \in \RRRRRR_{i,j} \text{ and } i+j \GEQSLANT \TEMPORARYI-1,\end{array}\right.\]
where 
\begin{align*}
    \textstyle \lambda_{k,\TEMPORARYI} & \textstyle = \frac{\binom{s-\TEMPORARYI+2}{\TEMPORARYI-1}a}{\binom{s-k+2}{\TEMPORARYI-1}p^{\TEMPORARYI-1}} \in \smash{\overline{\FF}}_p^\times,
    \\ \textstyle \lambda_{k,\TEMPORARYI,i,j} & \textstyle = \frac{(-1)^{\TEMPORARYI+i+j+1} (\TEMPORARYI-j-1) \binom{\TEMPORARYI-j-2}{i}\binom{s-\TEMPORARYI+j+2}{\TEMPORARYI-j-1}a}{(k-s-i(p-1)-2)p^{\TEMPORARYI-j-1}} \in \smash{\overline{\FF}}_p^\times.
\end{align*}
\end{theorem}

Thus our task is to prove theorems \ref{MAINTHM0a} and \ref{MAINTHM00a}.

\section{Combinatorics}
\label{secComb}

Throughout the proof we will refer to the combinatorial results in section 8 of~\cite{Ars1}. For convenience, we reproduce the statements here in the form we will use.

\begin{lemma}\label{part1X5}
Suppose throughout this lemma that 
\begin{align*}
    n,t,y \in \ZZ, \qquad
    b,d,k,l,w  \in \ZZ_{\GEQSLANT 0}, \qquad
    m,u,v  \in \ZZ_{\GEQSLANT 1}.
\end{align*}
\begin{enumerate}
    \item If \m{u \equiv v \bmod (p-1)p^{m-1}} then
\begin{align}
\label{ceqn1}\tag{\emph{c-a}}
\textstyle \TEMPORARYM_{u,n} \equiv \TEMPORARYM_{v,n} \bmod p^m.
\end{align}
\item Suppose that \m{u = t_u(p-1)+s_u} with \m{s_u = \MORE{u}}, so that \m{s_u \in \{1,\ldots,p-1\}} and \m{t_u \in \ZZ_{\GEQSLANT 0}}. Then \begin{align}
\label{ceqn2}\tag{\emph{c-b}}
\textstyle \TEMPORARYM_u = 1 + \VERIFYIF{u\equiv_{p-1} 0} + \frac {t_u}{s_u}p + \BIGO{t_up^2}.
\end{align}
\item  If \m{n\LEQSLANT 0} then
\begin{align}
\label{ceqn3}\tag{\emph{c-c}}
\textstyle \TEMPORARYM_{u,n} = %
\sum_{i=0}^{-n} (-1)^i \binom{-n}{i} \TEMPORARYM_{u-n-i,0}. %
\end{align} \item If \m{n\GEQSLANT 0} then
\begin{align}
\label{ceqn4}\tag{\emph{c-d}}\textstyle \TEMPORARYM_{u,n} \equiv (1+\VERIFYIF{u\equiv_{p-1}n\equiv_{p-1} 0})\binom{\MORE{u}}{\LESS{n}} %
 \bmod p.
\end{align}
\item If \m{u \GEQSLANT (b+l)d} and \m{l\GEQSLANT w} then
\begin{align}
\label{ceqn5}\tag{\emph{c-e}}\textstyle \sum_j (-1)^{j-b} \binom{l}{j-b}\binom{u-dj}{w} = \VERIFYIF{w=l}d^l.
\end{align}
\item If \m{\FORMX} is a formal variable  then
\begin{align}\label{ceqn8}\tag{\emph{c-f}}
   \textstyle  \binom{\FORMX}{t+l} \binom{t}{w} = \sum_v (-1)^{w-v} \binom{l+w-v-1}{w-v} \binom{\FORMX}{v} \binom{\FORMX-v}{t+l-v}.
\end{align}
Consequently, if \m{b+l\GEQSLANT d+w} then
\begin{align}
\label{ceqn6}\nonumber& \textstyle \sum_i %
 \binom{b-d+l}{i(p-1)+l}\binom{i(p-1)}{w} \\ & \textstyle \null\qquad = \sum_v (-1)^{w-v} \binom{l+w-v-1}{w-v} \binom{b-d+l}{v} \TEMPORARYM_{b-d+l-v,l-v}.\tag{\emph{c-g}}\end{align}
\item We have \begin{align}
\label{ceqn9}\tag{\emph{c-i}}
\textstyle \sum_{j} (-1)^j \binom{y}{j} \binom{y+l-j}{w-j} = (-1)^w \binom{w-l-1}{w}.
\end{align}
\item We have 
\begin{align}
\label{ceqn10}\tag{\emph{c-j}}
\textstyle \sum_{j} \binom{u-1}{j-1} \binom{-l}{j-w} = (-1)^{u-w} \binom{l-w}{u-w}.
\end{align}
\item We have
\begin{align}
\label{ceqn11}\tag{\emph{c-k}}
\textstyle \sum_{j} (-1)^{j} \binom{j}{b} \binom{l}{j-w} = (-1)^{l+w}\binom{w}{l+w-b}.
\end{align}
\end{enumerate}
\end{lemma}

\begin{lemma}\label{part2X5}
Let \m{\alpha \in \ZZ \cap [0,\ldots,\frac{r}{p+1}]} and let \m{\{D_i\}_{i\in \ZZ}} be a family of elements of \m{\ZZ_p} such that \m{D_i=0} for \m{i\not\in[0,\frac{r-\alpha}{p-1}]} and
\m{\textstyle \TEMPORARYD_{w} (D_\DARKCIRC) = 0} for all \m{0 \LEQSLANT w<\alpha}. Then \[\textstyle \sum_i D_i x^{i(p-1)+\alpha}y^{r-i(p-1)-\alpha} = \theta^\alpha h\] for some polynomial \m{h} with integer coefficients.
\end{lemma}

\begin{lemma}\label{part3X5} For \m{\alpha,\lambda,\mu \in \ZZ_{\GEQSLANT 0}} let
\[\textstyle L_\alpha(\lambda,\mu)\] be the 
\m{(\alpha+1)\times(\alpha+1)} matrix with %
 entries
    \[\textstyle L_{l,j} = \sum_{k=0}^{\alpha} \frac{j!}{l!} \left(\frac\mu\lambda\right)^k  \sfirst (l,k) \ssecond (k,j),\]
where \m{\sfirst (l,k)} are the Stirling numbers of the first kind and \m{\ssecond (k,j)} are the Stirling numbers of the second kind. Then
\[\textstyle L_\alpha(\lambda,\mu) \left(\binom{\lambda \FORMX}{0},\ldots,\binom{\lambda \FORMX}{\alpha}\right)^T = \left(\binom{\mu \FORMX}{0},\ldots,\binom{\mu \FORMX}{\alpha}\right)^T.\]
\end{lemma}

\begin{lemma}\label{lemmaComb2} 
For \m{\alpha \in \ZZ_{\GEQSLANT 0}} let \m{B_{\alpha}} be the 
\m{(\alpha+1)\times(\alpha+1)} matrix with %
 entries
    \[\textstyle B_{i,j} = j! \sum_{k,l=0}^{\alpha} \frac{(-1)^{i+l+k}}{l!} \binom{l}{i} (1-p)^{-k} \sfirst (l,k) \ssecond (k,j),\]
where \m{\sfirst (i,j)}  and \m{\ssecond (k,j)} are the Stirling numbers of the first and second kind, respectively. Let \m{\{X_{i,j}\}_{i,j\GEQSLANT 0}} be formal variables.
For \m{\beta \in \ZZ_{\GEQSLANT 0}} such that \m{\alpha\GEQSLANT \beta} let  \[\textstyle S({\alpha,\beta})=\left(S({\alpha,\beta})_{w,j}\right)_{0\LEQSLANT w,j\LEQSLANT\alpha}\] be the \m{(\alpha+1)\times(\alpha+1)} matrix with %
 entries
 \[\textstyle S({\alpha,\beta})_{w,j} = \sum_{i=1}^\beta X_{i,j} \binom{i(p-1)}{w}.\]
 Then \m{B_\alpha S({\alpha,\beta})} is zero outside the rows indexed \m{1,\ldots,\beta} and \[\textstyle (B_\alpha S({\alpha,\beta}))_{i,j} = X_{i,j}\]
for \m{i \in \{1,\ldots,\beta\}}.
\end{lemma}

\begin{lemma}\label{part2X3a}
For \m{u,v,c \in \ZZ} let us define
\[\textstyle F_{u,v,c}(X) = \sum_w (-1)^{w-c} \binom{w}{c} \binom{X}{w}^\partial \binom{X+u-w}{v-w} \in \QQ_p[X].\]
Then
\[\textstyle F_{u,v,c}(X) =  \binom{u}{v-c} \binom{X}{c}^\partial-\binom{u}{v-c}^\partial \binom{X}{c}.\]
\end{lemma}

\begin{lemma}\label{part4X5} Let \m{\FORMX} and \m{\FORMY} denote formal variables, and let
\[\textstyle c_j = (-1)^j \alpha! \left( \frac{\FORMX+j+1}{j+1}\binom{\FORMY}{\alpha-j-1} + \binom{\FORMY}{\alpha-j}\right)\in\QQ[\FORMX,\FORMY]\subset \QQ(\FORMX,\FORMY)\]
be polynomials over \m{\QQ} of degrees \m{\alpha-j}, for \m{1\LEQSLANT j\LEQSLANT \alpha}. Let \[\textstyle M = (M_{w,j})_{0\LEQSLANT w,j\LEQSLANT \alpha}\] be the \m{(\alpha+1)\times(\alpha+1)} matrix over \m{\QQ(\FORMX,\FORMY)} with entries
 \begin{align*}
        M_{w,0} & =\textstyle (-1)^w\frac{(\FORMY-\FORMX)\FORMX_w}{\FORMY_{w+1}}, %
        \\ M_{w,j} & = \textstyle \sum_v (-1)^{w-v} \binom{j+w-v-1}{w-v} \binom{\FORMX+j}{v} \left(\binom{\FORMY+j-v}{j-v}-\binom{\FORMX+j-v}{j-v}\right), %
\end{align*}
for \m{0\LEQSLANT w\LEQSLANT \alpha} and \m{0<j\LEQSLANT\alpha}. Then the first \m{\alpha-1} entries of
\[\textstyle Mc  = M(\FORMY_{\alpha},c_1,\ldots,c_\alpha)^T = (d_0,\ldots,d_\alpha)^T\] are zero, and \m{d_\alpha=\frac{(\FORMY-\FORMX)_{\alpha+1}}{\FORMY-\alpha}}. \end{lemma}

\begin{lemma}\label{lemmaComb4} 
Suppose that \m{s,\alpha,\beta\in\ZZ} are such that
\[\textstyle 1 \LEQSLANT \beta \LEQSLANT \alpha \LEQSLANT \frac{s}{2}-2 \LEQSLANT \frac{p-5}{2}.\]
Let \m{B=B_\alpha} denote the matrix defined in lemma~\ref{lemmaComb2}.
Let \m{M} denote the \m{(\alpha+1)\times(\alpha+1)} matrix with entries in \m{\FF_p} such that if \m{i\in\{1,\ldots,\beta\}} and \m{j\in\{0,\ldots,\alpha\}} then
\[\textstyle M_{i,j} = \binom{\beta}{i} \cdot \left\{\begin{array}{rl}
    \binom{s-\alpha-\beta+i}{i}^{-1} (-1)^{i+1} & \text{if } j = 0, \\
    \binom{s-\alpha-\beta+j}{j-i} & \text{if }j>0,
\end{array}\right.\]
and if \m{i\in\{0,\ldots,\alpha\}\backslash\{1,\ldots,\beta\}} and \m{j\in\{0,\ldots,\alpha\}} then \m{M_{i,j}} is the reduction modulo \m{p} of
 \begin{align*}
       & \textstyle p^{-\VERIFYIF{j=0}}  \sum_{w=0}^\alpha B_{i,w} \sum_v (-1)^{w-v}  \binom{j+w-v-1}{w-v} \binom{s+\beta(p-1)-\alpha+j}{v}^\partial  
       \\
        & \textstyle \null\qquad\qquad\qquad\qquad\qquad\qquad\qquad\qquad\qquad\qquad \cdot\sum_{u=0}^{\beta} \binom{s+\beta(p-1)-\alpha+j-v}{u(p-1)+j-v}
        \\
        & \textstyle \null\qquad - \VERIFYIF{i=0} p^{-\VERIFYIF{j=0}} \textstyle  \binom{s+\beta(p-1)-\alpha+j} {j}^\partial \\
        & \textstyle \null\qquad  -\VERIFYIF{j=0} \sum_{w=0}^\alpha B_{i,w} (-1)^{w}\binom{s+\beta(p-1)-\alpha}{w}  \frac{w!}{(s-\alpha)_{w+1}}.
\end{align*}
Then there is a solution of
\[M (z_0,\ldots,z_\alpha)^T = (1,0,\ldots,0)^T\] such that \m{z_0 \ne 0}. 
\end{lemma}

Now let us prove some additional combinatorial results.

\begin{lemma}\label{lemmaComb555}
Suppose that \m{\alpha\in\ZZ_{\GEQSLANT0}}. For \m{w,j \in \{0,\ldots,\alpha\}} let \[\textstyle F_{w,j}(\polyr,\polys) \in \FF_p[\polyr,\polys]\] denote the polynomial
\[\textstyle \sum_v (-1)^{w-v} \binom{j+w-v-1}{w-v}\binom{\polyr-\alpha+j}{v} \left(\binom{\polys-\alpha+j-v}{j-v}-\binom{\polyr-\alpha+j-v}{j-v}\right).\]
Note that this depends on \m{\alpha}. Then
\begin{align*} \textstyle
     \sum_{j=1}^\alpha (-1)^{\alpha-j} \binom{\polys-\alpha+1}{\alpha-j} F_{w,j}(\polyr,\polys)  = 
  (-1)^\alpha(\VERIFYIF{w=\alpha}-\VERIFYIF{w=0}) \binom{\polys-\polyr}{\alpha}.
\end{align*}
\end{lemma}

\begin{proof}{Proof. }
   Both sides of the equation we want to prove have degree \m{\alpha} and the coefficient of \m{\polyr^\alpha} on each side is \m{\frac{1}{\alpha!}\left(\VERIFYIF{w=\alpha}-\VERIFYIF{w=0}\right)}. So the two sides are equal if they are equal when evaluated at the points \m{(\polyr,\polys)} such that
\[\textstyle (\polyr,\polys) \in \left\{ (u+\gamma(p-1)+\alpha,u+\alpha) \,|\, u \in \{0,\ldots,\alpha\} , \, \gamma \in \{0,\ldots,\alpha-1\}\right\}.\]
    The right side is zero when evaluated at these points, and
    \[\textstyle F_{w,j} (u+\gamma(p-1)+\alpha,u+\alpha) = \sum_{i=1}^\gamma \binom{u+\gamma(p-1)+j}{i(p-1)+j} \binom{i(p-1)}{w}\]
    by~\eqref{ceqn6}. Thus we want to show that
    \[\textstyle \sum_{j=1}^\alpha (-1)^{\alpha-j} \binom{u+1}{\alpha-j} \sum_{i=1}^\gamma \binom{u+\gamma(p-1)+j}{i(p-1)+j}\binom{i(p-1)}{w} = \BIGO{p}\] for \m{0\LEQSLANT u,w\LEQSLANT \alpha} and \m{0\LEQSLANT \gamma< \alpha}. Since \[\textstyle\binom{u+\gamma(p-1)+j}{i(p-1)+j}\binom{i(p-1)}{w} = \binom{\gamma}{i}\binom{u+j-\gamma}{j-i}\binom{-i}{w}+\BIGO{p},\] that is equivalent to
\[\textstyle \sum_{i,j>0} (-1)^{\alpha+w-i} \binom{u+1}{\alpha-j} \binom{\gamma}{i}\binom{\gamma-u-i-1}{j-i}\binom{i+w-1}{w} = \BIGO{p}.\] 
This follows from the facts that
\[\textstyle \sum_{j>0} \binom{u+1}{\alpha-j}\binom{\gamma-u-i-1}{j-i} = \binom{\gamma-i}{\alpha-i}\]
for \m{i>0} by Vandermonde's convolution formula, and \[\textstyle  \binom{\gamma}{i} \binom{\gamma-i}{\alpha-i} = \binom{\alpha}{i} \binom{\gamma}{\alpha}=0\] since \m{\gamma\in\{0,\ldots,\alpha-1\}}.
\end{proof}

  \section{Computing  \texorpdfstring{\m{\overline\Theta_{k,a}}}{\textTheta}}
    \label{secComp4}
    
    Throughout the proof we use the results from section~9 of~\cite{Ars1}, which we reproduce here without proofs for convenience.
    
    \begin{lemma}\label{LEMMAX2}
{
Suppose that \m{\alpha \in \{0,\ldots,\TEMPORARYI-1\}}.
\begin{enumerate}
    \item \label{111}
We have
  \begin{align*}
      \textstyle & \textstyle (T-a) \left(1 \sqbrackets{}{KZ}{\QQp}{ \theta^{n}x^{\alpha-n}y^{r-np-\alpha}}\right)
      \\ & \textstyle \null\qquad 
      = \sum_j (-1)^j \binom{n}{j} p^{j(p-1)+\alpha} \begin{Lmatrix}
      1 & 0 \\ 0 & p
      \end{Lmatrix} \sqbrackets{}{KZ}{\QQp}{x^{j(p-1)+\alpha}y^{r-j(p-1)-\alpha}} 
      \\ & \textstyle \null\qquad\qquad - a \sum_j (-1)^j \binom{n-\alpha}{j} \sqbrackets{}{KZ}{\QQp}{\theta^\alpha x^{j(p-1)}y^{r-j(p-1)-\alpha(p+1)}} + \BIGO{p^n}.
  \end{align*}
\item \label{222} The submodule \m{\Img(T-a)\subset\ind_{KZ}^G \TILDE\SIGMA_{r}} contains
  \begin{align*}
      & \textstyle \sum_{i} \left(\sum_{l=\beta-\gamma}^{\beta} C_l \binom{r-\beta+l}{i(p-1)+l}\right)\sqbrackets{}{KZ}{\QQp}{x^{i(p-1)+\beta}y^{r-i(p-1)-\beta}} \\ & \textstyle \null\qquad\qquad\qquad\qquad\qquad\qquad\qquad\qquad\qquad+ \BIGO{ap^{-\beta+v_C}+p^{p-1}}
  \end{align*}
  for all \m{0 \LEQSLANT \beta \LEQSLANT \gamma < \TEMPORARYI} and all families \m{\{C_l\}_{l\in \ZZ}} of elements of \m{\ZZ_p}, where
  \[\textstyle v_C = \min_{\beta-\gamma\LEQSLANT l \LEQSLANT \beta}(\PADICVALUATION(C_l)+l).\]
  The \m{\BIGO{ap^{-\beta+v_C}+p^{p-1}}} term is equal to  \m{\BIGO{p^{p-1}}} plus
  \[\textstyle -\frac{ap^{-\beta}}{p-1} \sum_{l=\beta-\gamma}^{\beta} C_{l}p^l\sum_{0\ne\MU\in\FF_p} [\MU]^{-l}\begin{Lmatrix}
p & [\MU] \\ 0 & 1
\end{Lmatrix}  \sqbrackets{}{KZ}{\QQp}{\theta^n  x^{\beta-l-n}y^{r-np-\beta+l}}.\]
\end{enumerate}
}
\end{lemma}

\begin{lemma}\label{LEMMAX3}
{
Suppose that \m{\alpha \in \ZZ} and \m{v \in \QQ} are such that
\begin{align*}\textstyle \alpha & \textstyle \in\{0,\ldots,\TEMPORARYI-1\}, \\ \textstyle v & \textstyle \LEQSLANT \PADICVALUATION(\TEMPORARYD_{\alpha}(D_\DARKCIRC)),\\v' \defeq \min\{\PADICVALUATION(a)-\alpha,v\} & \textstyle  \LEQSLANT \PADICVALUATION(\TEMPORARYD_{w}(D_\DARKCIRC)) \text{ for } \alpha < w < 2\TEMPORARYI-\alpha,
\\ v' & \textstyle <\PADICVALUATION(\TEMPORARYD_{w}(D_\DARKCIRC)) \text{ for } 0\LEQSLANT w<\alpha.
\end{align*} If, for \m{j \in \ZZ},
\[\textstyle\DELTAA_j \defeq (-1)^{j-1} (1-p)^{-\alpha} \binom{\alpha}{j-1} \TEMPORARYD_{\alpha}(D_\DARKCIRC),\]
then \m{v \LEQSLANT\PADICVALUATION(\TEMPORARYD_{\alpha}(\DELTAA_\DARKCIRC))\LEQSLANT \PADICVALUATION(\DELTAA_j)} for all \m{j \in \ZZ}, and
\begin{align*}
      \textstyle & \textstyle \sum_i (\DELTAA_i - D_i) \sqbrackets{}{KZ}{\QQp}{x^{i(p-1)+\alpha}y^{r-i(p-1)-\alpha}}
      \\ & \textstyle \null\qquad 
      =  \VERIFYIF{\alpha\LEQSLANT s}(-1)^{n+1} D_{\frac{r-s}{p-1}}
 \sqbrackets{}{KZ}{\QQp}{\theta^n x^{r-np-s+\alpha}y^{s-\alpha-n}} 
            \\ & \textstyle \null\qquad \qquad - D_0
 \sqbrackets{}{KZ}{\QQp}{\theta^n x^{\alpha-n}y^{r-np-\alpha}}  %
  \\ & \textstyle \null\qquad \qquad \qquad + E \sqbrackets{}{KZ}{\QQp}{\theta^{\alpha+1}h} + %
  F \sqbrackets{}{KZ}{\QQp}{h'} + \mathrm{ERR}_1 + \mathrm{ERR}_2,
  \end{align*}
  for some \m{\mathrm{ERR}_1} and \m{\mathrm{ERR}_2} such that
  \[\textstyle \mathrm{ERR}_1 \in\Img(T-a) \text{ and }\mathrm{ERR}_2 = \BIGO{p^{\TEMPORARYI-\PADICVALUATION(a)+v}+p^{\TEMPORARYI-\alpha}},\]
 some polynomials %
 \m{h} and \m{h'}, and some %
  \m{E,F\in\QQp} such that %
  \m{\PADICVALUATION(E)\GEQSLANT v'} and \m{\PADICVALUATION(F)>v'}.
}
\end{lemma}

\begin{lemma}\label{LEMMAX4}
{
Let  \m{\{C_l\}_{l\in \ZZ}} be any family of elements of \m{\ZZ_p}. Suppose that 
\m{\alpha \in \{0,\ldots, \TEMPORARYI-1\}} and \m{v \in \QQ}, and suppose that the constants
 \[\textstyle D_i \defeq %
 \VERIFYIF{i=0} C_{-1} + \VERIFYIF{0<i(p-1)<r-2\alpha}\sum_{l=0}^\alpha C_{l} \binom{r-\alpha+l}{i(p-1)+l}\]
satisfy the conditions of lemma~\ref{LEMMAX3}, \IE 
\begin{align*}\textstyle  \textstyle v & \textstyle \LEQSLANT \PADICVALUATION(\TEMPORARYD_{\alpha}(D_\DARKCIRC)),\\v' \defeq \min\{\PADICVALUATION(a)-\alpha,v\} & \textstyle  \LEQSLANT \PADICVALUATION(\TEMPORARYD_{w}(D_\DARKCIRC)) \text{ for } \alpha < w < 2\TEMPORARYI-\alpha,
\\ v' & \textstyle <\PADICVALUATION(\TEMPORARYD_{w}(D_\DARKCIRC)) \text{ for } 0\LEQSLANT w<\alpha.
\end{align*}
Moreover, suppose that \m{C_0} is a unit. Let
\[\textstyle\TEMPORARYD' \defeq (1-p)^{-\alpha}\TEMPORARYD_{\alpha}(D_\DARKCIRC) - C_{-1}   .\]
Suppose that \m{\PADICVALUATION(C_{-1})\GEQSLANT \PADICVALUATION(\TEMPORARYD')}.
\begin{enumerate}
    \item \label{7111} If \m{\PADICVALUATION(\TEMPORARYD') \LEQSLANT v'}
then there is some element  \m{\mathrm{gen}_1 \in \SCR I_a} that
represents a generator of \m{\FILT N_\alpha}. \\ \item \label{7222} If \m{\PADICVALUATION(a)-\alpha < v} then there is some element \m{\mathrm{gen}_2 \in \SCR I_a}  that represents a generator of a finite-codimensional submodule of 
\[\textstyle T\left(\ind_{KZ}^G\QUOTIENTI(\alpha)\right) =T\left(\FILT N_\alpha/ \ind_{KZ}^G\SUBMODULE(\alpha)\right),\]
 where \m{T} denotes the  endomorphism  of \m{\ind_{KZ}^G\QUOTIENTI(\alpha)} corresponding to the double coset of \m{\begin{Smatrix}
p & 0 \\ 0 & 1
\end{Smatrix}}.
\end{enumerate}
}
\end{lemma}

Let us now prove the following additional results.

    \begin{lemma}\label{LEMMATHM}
Let \m{\{C_l\}_{l\in \ZZ}} be any family of elements of~\m{\ZZ_p}. Suppose that \m{\alpha \in \ZZ} and \m{v \in \QQ} and the constants
 \[\textstyle D_i \defeq \VERIFYIF{i=0} C_{-1} + \VERIFYIF{0<i(p-1)<r-2\alpha}\sum_{l=0}^\alpha C_{l} \binom{r-\alpha+l}{i(p-1)+l}\] are such that
\begin{align*}\textstyle \alpha & \textstyle \in\{0,\ldots,\TEMPORARYI-1\}, \\ \textstyle v & \textstyle \LEQSLANT \PADICVALUATION(\TEMPORARYD_{\alpha}(D_\DARKCIRC)),
\\ v' \defeq \min\{\PADICVALUATION(a)-\alpha,v\}  & \textstyle <\PADICVALUATION(\TEMPORARYD_{\TEMPORARYAAAA}(D_\DARKCIRC)) \text{ for } 0\LEQSLANT \TEMPORARYAAAA<\alpha.
\end{align*} 
Let
\[\textstyle\TEMPORARYD' \defeq (1-p)^{-\alpha} \TEMPORARYD_{\alpha}(D_\DARKCIRC) - C_{-1}.\]
Then \m{\Img(T-a)} contains
\begin{align}& \textstyle (\TEMPORARYD'+C_{-1}) \sqbrackets{}{KZ}{\QQp}{\theta^\alpha x^{p-1}y^{r-\alpha(p+1)-p+1}}  + C_{-1}  \sqbrackets{}{KZ}{\QQp}{\theta^n x^{\alpha-n}y^{r-np-\alpha}}\nonumber
\\ & \textstyle\qquad \label{EQIMPX} + \sum_{\TEMPORARYA=\alpha+1}^{2\TEMPORARYI-\alpha-1} E_\TEMPORARYA 
 \sqbrackets{}{KZ}{\QQp}{\theta^{\TEMPORARYA}h_\TEMPORARYA %
 }  +  F \sqbrackets{}{KZ}{\QQp}{h'}
+ H%
,\end{align}
for some \m{h_\TEMPORARYA%
,h',E_\TEMPORARYA
,F,H} such that 
\begin{enumerate}
    \item\label{mainlemp1} \m{E_\TEMPORARYA = \TEMPORARYD_{\TEMPORARYA}(D_\DARKCIRC) + \BIGO{p^v} \cup \BIGO{\TEMPORARYD_{\alpha+1}(D_\DARKCIRC)} \cup\cdots\cup \BIGO{\TEMPORARYD_{\TEMPORARYA-1}(D_\DARKCIRC)}}, \\
     \item\label{mainlemp1a} if \m{\LESS{\TEMPORARYA+\alpha-s}\LEQSLANT\LESS{2\TEMPORARYA-s}\ne 0} then the reduction modulo \m{\maxideal} of \m{\textstyle \theta^{\TEMPORARYA}h_\TEMPORARYA} generates \m{N_\TEMPORARYA},\\
    \item\label{mainlemp2} \m{\PADICVALUATION(F)>v'}, and\\
    \item\label{mainlemp3} \m{H = \BIGO{p^{\TEMPORARYI-\PADICVALUATION(a)+v}+p^{\TEMPORARYI-\alpha}}} and if \m{\PADICVALUATION(a) -\alpha< v } then \[\textstyle \frac{1-p}{ap^{-\alpha}} H = g  \sqbrackets{}{KZ}{\QQp}{\theta^n x^{\alpha-n}y^{r-np-\alpha}} + \BIGO{p^{\TEMPORARYI-\PADICVALUATION(a)}}\]
    with
    \[\textstyle g = \sum_{\lambda\in\FF_p} C_0 \begin{Lmatrix}
    p & [\lambda] \\ 0 & 1\end{Lmatrix} + A \begin{Lmatrix}
    p & 0 \\ 0 & 1\end{Lmatrix} + \VERIFYIF{r\equiv_{p-1} 2\alpha} B  \begin{Lmatrix}
    0 & 1 \\ p & 0 \end{Lmatrix},\]
    where \[\textstyle A = -C_{-1}+\sum_{l=1}^\alpha C_l \binom{r-\alpha+l}{l}\] and \[\textstyle B=\sum_{l=0}^\alpha C_l \binom{r-\alpha+l}{s-\alpha}.\]
\end{enumerate}
\end{lemma}

\begin{proof}{Proof. } This lemma is essentially shown under a stronger hypothesis as lemma~\ref{LEMMAX4}. The stronger hypothesis consists of the three extra conditions that \m{\PADICVALUATION(\TEMPORARYD_{\TEMPORARYAAAA}(D_\DARKCIRC))\GEQSLANT \min\{\PADICVALUATION(a)-\alpha,v\}} for all \m{\alpha<\TEMPORARYAAAA<2\TEMPORARYI-\alpha}, that
\m{C_0\in \ZZ_p^\times}, and that \m{\PADICVALUATION(C_{-1})\GEQSLANT \PADICVALUATION(\TEMPORARYD')}.
These extra conditions are not used in the actual construction of the element in~\eqref{EQIMPX}, rather they are there to ensure that
\m{\PADICVALUATION(E_\TEMPORARYA) \GEQSLANT \min\{\PADICVALUATION(a)-\alpha,v\}} for all \m{\alpha<\TEMPORARYA<2\TEMPORARYI-\alpha}, that  the coefficient of \m{\begin{Smatrix}p & [\lambda] \\ 0 & 1\end{Smatrix}} in \m{g} is invertible, and that we get an integral element once we divide the element
\[\textstyle (\TEMPORARYD'+C_{-1}) \sqbrackets{}{KZ}{\QQp}{\theta^\alpha x^{p-1}y^{r-\alpha(p+1)-p+1}}  + C_{-1}  \sqbrackets{}{KZ}{\QQp}{\theta^n x^{\alpha-n}y^{r-np-\alpha}}\]
by \m{\TEMPORARYD'}.
Therefore we still get the existence of the element in~\eqref{EQIMPX} without these extra conditions, and to complete the proof of lemma~\ref{LEMMATHM} we need to verify the properties of \m{h_\TEMPORARYA,E_\TEMPORARYA,F,H,A}, and \m{B} claimed in~(\ref{mainlemp1}), (\ref{mainlemp1a}), (\ref{mainlemp2}), and (\ref{mainlemp3}).
 The \m{h_\TEMPORARYA} and \m{E_\TEMPORARYA} come from the proof of lemma~\ref{LEMMAX3}, and \m{E_\TEMPORARYA 
 \sqbrackets{}{KZ}{\QQp}{\theta^{\TEMPORARYA}h_\TEMPORARYA}} is
 \[\textstyle X_\TEMPORARYA 
 \sqbrackets{}{KZ}{\QQp}{\sum_{\lambda\ne 0} [-\lambda]^{r-\alpha-\TEMPORARYA} \begin{Lmatrix}1 & [\lambda] \\ 0 & 1\end{Lmatrix} (-\theta)^{n} x^{r-np-\TEMPORARYA}y^{\TEMPORARYA-n}},\]
with the notation for \m{X_\TEMPORARYA} from the proof of lemma~\ref{LEMMAX3}. Let \m{E_\TEMPORARYA=(-1)^{\TEMPORARYA+1}X_\TEMPORARYA}. Then condition~(\ref{mainlemp1}) is satisfied directly from the definition of \m{X_\TEMPORARYA}. Let
 \[\textstyle h_\TEMPORARYA = (-1)^{\TEMPORARYA+1} \sum_{\lambda\ne 0} [-\lambda]^{r-\alpha-\TEMPORARYA} \begin{Lmatrix}1 & [\lambda] \\ 0 & 1\end{Lmatrix} (-\theta)^n x^{r-np-\TEMPORARYA}y^{\TEMPORARYA-n}.\] This reduces modulo \m{\maxideal} to the element
\[\textstyle (-1)^{\TEMPORARYA} \sum_{\lambda\ne 0} [-\lambda]^{r-\alpha-\TEMPORARYA} \begin{Lmatrix}1 & [\lambda] \\ 0 & 1\end{Lmatrix} {Y^{\LESS{2\TEMPORARYA-r}}} = (-1)^{s-\alpha+1} \binom{\LESS{2\TEMPORARYA-s}}{\LESS{\TEMPORARYA+\alpha-s}} X^{\LESS{\TEMPORARYA+\alpha-s}} Y^{\LESS{\TEMPORARYA-\alpha}} \]
of  \[\textstyle \sigma_{\LESS{2\TEMPORARYA-r}}(r-\TEMPORARYA) \cong I_{r-2\TEMPORARYA}(\TEMPORARYA) / \sigma_{\MORE{r-2\TEMPORARYA}}(\TEMPORARYA)=\QUOTIENTI(\TEMPORARYA).\] This element is non-trivial and generates \m{N_\TEMPORARYA} if \m{\LESS{\TEMPORARYA+\alpha-s} \LEQSLANT \LESS{2\TEMPORARYA-s} \ne 0}, since then %
 \m{X^{\LESS{\TEMPORARYA+\alpha-s}} Y^{\LESS{\TEMPORARYA-\alpha}}} generates \m{N_\TEMPORARYA}. This verifies condition~(\ref{mainlemp1a}). Condition~(\ref{mainlemp2}) follows from the assumption \m{v'<\PADICVALUATION(\TEMPORARYD_{w}(D_\DARKCIRC))} for \m{0\LEQSLANT w < \alpha}, as in the proof of lemma~\ref{LEMMAX3}. Finally, condition~(\ref{mainlemp3}) follows from the description of the error term in lemma~\ref{LEMMAX2}, as in the proof of lemma~\ref{LEMMAX4}.
\end{proof}

\begin{corollary}\label{CORO1}
Let \m{\{C_l\}_{l\in \ZZ}} be any family of elements of~\m{\ZZ_p}. Suppose that \m{\alpha \in \{0,\ldots,\TEMPORARYI-1\}} and \m{v \in \QQ}, and suppose that the constants
 \[\textstyle D_i \defeq %
 \VERIFYIF{i=0} C_{-1} + \VERIFYIF{0<i(p-1)<r-2\alpha}\sum_{l=0}^\alpha C_{l} \binom{r-\alpha+l}{i(p-1)+l}\] are such that
\begin{align*}\textstyle  \textstyle v & \textstyle \LEQSLANT \PADICVALUATION(\TEMPORARYD_{\alpha}(D_\DARKCIRC)),
\\
v' \defeq \min\{\PADICVALUATION(a)-\alpha,v\}  & \textstyle \LEQSLANT \PADICVALUATION(\TEMPORARYD_{\TEMPORARYAAAA}(D_\DARKCIRC)) \text{ for } \alpha<\TEMPORARYAAAA<2\TEMPORARYI-\alpha,
\\ v'  & \textstyle <\PADICVALUATION(\TEMPORARYD_{\TEMPORARYAAAA}(D_\DARKCIRC)) \text{ for } 0\LEQSLANT w<\alpha.
\end{align*} 
Suppose also that \m{\PADICVALUATION(a)\not\in \ZZ}. Let \begin{align*} \textstyle\TEMPORARYD' &\textstyle\defeq (1-p)^{-\alpha} \TEMPORARYD_{\alpha}(D_\DARKCIRC) - C_{-1}, \\ \textstyle\veeC & \textstyle\defeq-C_{-1}+\sum_{l=1}^\alpha C_l \binom{r-\alpha+l}{l}.\end{align*}
 If \m{\star} then \m{\ast} is trivial modulo \m{\SCR I_a}, for each of the following pairs \[\textstyle(\star,\ast) = (\text{condition}, \text{representation}).\]

\begin{enumerate}
    \item\label{CORO1.1} \m{\Big(\PADICVALUATION(\TEMPORARYD') \LEQSLANT \min\{\PADICVALUATION(C_{-1}),v'\},\,\FILT N_\alpha\Big)}.\NEWLINEA
        \item\label{CORO1.2}  \m{\Big(v=\PADICVALUATION(C_{-1}) < \min\{\PADICVALUATION(\TEMPORARYD'),\PADICVALUATION(a)-\alpha\},\,\ind_{KZ}^G \SUBMODULE(\alpha)\Big)}.\NEWLINEA
            \item\label{CORO1.3}  \m{\Big(\PADICVALUATION(a)-\alpha < v\LEQSLANT \PADICVALUATION(C_{-1}) \annd \textstyle \veeC\in \ZZ_p^\times \annd C_0 \not\in\ZZ_p^\times \annd \LESS{2\alpha-r}> 0,\, \FILT N_\alpha\Big)}.\NEWLINEA
                                \item\label{CORO1.4} \m{\Big(\PADICVALUATION(a)-\alpha < v\LEQSLANT \PADICVALUATION(C_{-1})  \annd \textstyle %
   \veeC \in \ZZ_p^\times,\,\ind_{KZ}^G \QUOTIENTI(\alpha)\Big)}.\NEWLINEA
                \item\label{CORO1.5} \m{\Big(\PADICVALUATION(a)-\alpha < v\LEQSLANT \PADICVALUATION(C_{-1}) \annd C_0 \in \ZZ_p^\times,\,\mathbf{r}_1\Big)}, where \[\textstyle\mathbf{r}_1\] is a finite-codimensional submodule of \[\textstyle T(\ind_{KZ}^G \QUOTIENTI(\alpha)).\]
\end{enumerate}
\end{corollary}

\begin{proof}{Proof. } %
  There is one extra condition imposed in addition to the conditions from lemma~\ref{LEMMATHM}: that
  \[\textstyle v' \defeq \min\{\PADICVALUATION(a)-\alpha,v\}  \textstyle \LEQSLANT \PADICVALUATION(\TEMPORARYD_{\TEMPORARYAAAA}(D_\DARKCIRC)) \text{ for } \alpha<\TEMPORARYAAAA<2\TEMPORARYI-\alpha,\] and it ensures that
\m{\PADICVALUATION(E_\TEMPORARYA) \GEQSLANT v'} for all \m{\alpha<\TEMPORARYA<2\TEMPORARYI-\alpha}. Lemma~\ref{LEMMATHM} implies that the element in~\eqref{EQIMPX} is in \m{\Img(T-a)}. Let us call this element \m{\gamma}.

(\ref{CORO1.1}) 
The condition \m{\PADICVALUATION(\TEMPORARYD') \LEQSLANT \min\{\PADICVALUATION(C_{-1}),v'\}} ensures that if we divide \m{\gamma} by \m{\TEMPORARYD'} then the resulting element reduces modulo~\m{\maxideal} to a representative of a generator of \m{\FILT N_\alpha}.

(\ref{CORO1.2}) 
The condition \m{v=\PADICVALUATION(C_{-1}) < \min\{\PADICVALUATION(\TEMPORARYD'),\PADICVALUATION(a)-\alpha\}} ensures that if we divide \m{\gamma} by \m{C_{-1}} then the resulting element   reduces modulo~\m{\maxideal} to a representative of a generator of \m{\ind_{KZ}^G \SUBMODULE(\alpha)}.

(\ref{CORO1.3}, \ref{CORO1.4}, \ref{CORO1.5}) 
The condition \m{\PADICVALUATION(a)-\alpha < v\LEQSLANT \PADICVALUATION(C_{-1}) } ensures that the term with the dominant valuation in~\eqref{EQIMPX} is \m{H}, so we can divide \m{\gamma} by \m{ap^{-\alpha}} and obtain the element \m{L + \BIGO{p^{\TEMPORARYI-\PADICVALUATION(a)}}}, where \m{L} is defined by
\[\textstyle L \defeq  \left(\sum_{\LAMBDA\in\FF_p}C_0\begin{Lmatrix}
p & [\LAMBDA] \\ 0 & 1
\end{Lmatrix} + A \begin{Lmatrix}
p & 0 \\0 & 1
\end{Lmatrix} +\VERIFYIF{r\equiv_{p-1}2\alpha} B \begin{Lmatrix}
0 & 1 \\ p& 0
\end{Lmatrix}\right) \sqbrackets{}{KZ}{\QQp}{\theta^n x^{\alpha-n}y^{r-np-\alpha}}\] with \m{A} and \m{B} as  in lemma~\ref{LEMMATHM}. This element \m{L} is in \m{\Img(T-a)}, and it
 reduces modulo~\m{\maxideal} to a representative of
\[\textstyle \left(\sum_{\LAMBDA\in\FF_p}C_0\begin{Lmatrix}
p & [\LAMBDA] \\ 0 & 1
\end{Lmatrix} + A \begin{Lmatrix}
p & 0 \\0 & 1
\end{Lmatrix} +\VERIFYIF{r\equiv_{p-1}2\alpha} (-1)^{r-\alpha} B \begin{Lmatrix}
1 & 0 \\0 & p
\end{Lmatrix}\right)\sqbrackets{}{KZ}{\FFp}{X^{\LESS{2\alpha-r}}}.\] %
 As shown in the proof of lemma~\ref{LEMMAX4}, if \m{C_0\in \ZZ_p^\times} then this element always generates a finite-codimensional submodule of \[\textstyle T(\ind_{KZ}^G \QUOTIENTI(\alpha)),\] and if additionally 
 \m{A\ne 0} (over \m{\FF_p}) then in fact we have the stronger conclusion that it generates \[\textstyle \ind_{KZ}^G \QUOTIENTI(\alpha).\] Suppose on the other hand that \m{C_0=\BIGO{p}} and \m{A\in\ZZ_p^\times}. In that case we assume that \m{\LESS{2\alpha-r}>0} and therefore the reduction modulo \m{\maxideal} of \m{L} represents a generator of \m{\FILT N_\alpha}. %
\end{proof}\begin{corollary}\label{CORO2}
Let \m{\{C_l\}_{l\in \ZZ}} be any family of elements of~\m{\ZZ_p}. Suppose that \m{\alpha \in \{0,\ldots,\TEMPORARYI-1\}} and \m{v \in \QQ}, and suppose that the constants
 \[\textstyle D_i \defeq %
 \VERIFYIF{i=0} C_{-1} + \VERIFYIF{0<i(p-1)<r-2\alpha}\sum_{l=0}^\alpha C_{l} \binom{r-\alpha+l}{i(p-1)+l}\] are such that
\begin{align*}\textstyle  \textstyle v & \textstyle \LEQSLANT \PADICVALUATION(\TEMPORARYD_{\alpha}(D_\DARKCIRC)),
\\
v' \defeq \min\{\PADICVALUATION(a)-\alpha,v\}  & \textstyle \LEQSLANT \PADICVALUATION(\TEMPORARYD_{\TEMPORARYAAAA}(D_\DARKCIRC)) \text{ for } \alpha<\TEMPORARYAAAA<2\TEMPORARYI-\alpha,
\\ v'  & \textstyle <\PADICVALUATION(\TEMPORARYD_{\TEMPORARYAAAA}(D_\DARKCIRC)) \text{ for } 0\LEQSLANT w<\alpha.
\end{align*} 
Suppose also that \m{\PADICVALUATION(a)\in \ZZ}. Let \begin{align*} \textstyle\TEMPORARYD' &\textstyle\defeq (1-p)^{-\alpha} \TEMPORARYD_{\alpha}(D_\DARKCIRC) - C_{-1}, \\ \textstyle\veeC & \textstyle\defeq-C_{-1}+\sum_{l=1}^\alpha C_l \binom{r-\alpha+l}{l}.\end{align*} If \m{\star} then \m{\ast} is trivial modulo \m{\SCR I_a}, for each of the following pairs \[\textstyle(\star,\ast) = (\text{condition}, \text{representation}).\]

\begin{enumerate}
    \item\label{CORO1.1b} \m{\Big(\PADICVALUATION(\TEMPORARYD') \LEQSLANT \min\{\PADICVALUATION(C_{-1}),v\} \annd \PADICVALUATION(\TEMPORARYD')  < \PADICVALUATION(a)-\alpha,\,\FILT N_\alpha\Big)}.\NEWLINEA
        \item\label{CORO1.2b}  \m{\Big(v=\PADICVALUATION(C_{-1}) < \min\{\PADICVALUATION(\TEMPORARYD'),\PADICVALUATION(a)-\alpha\},\,\ind_{KZ}^G \SUBMODULE(\alpha)\Big)}.\NEWLINEA
            \item\label{CORO1.3b}  \m{\Big(\PADICVALUATION(a)-\alpha < v\LEQSLANT \PADICVALUATION(C_{-1}) \annd \textstyle \veeC\in \ZZ_p^\times \annd C_0 \not\in\ZZ_p^\times \annd \LESS{2\alpha-r}> 0,\, \FILT N_\alpha\Big)}.\NEWLINEA
                                \item\label{CORO1.4b} \m{\Big(\PADICVALUATION(a)-\alpha < v\LEQSLANT \PADICVALUATION(C_{-1})  \annd \textstyle %
   \veeC \in \ZZ_p^\times,\,\ind_{KZ}^G \QUOTIENTI(\alpha)\Big)}.\NEWLINEA
                \item\label{CORO1.5b} \m{\Big(\PADICVALUATION(a)-\alpha < v\LEQSLANT \PADICVALUATION(C_{-1}) \annd C_0 \in \ZZ_p^\times,\,\mathbf{r}_2\Big)},  where \[\textstyle\mathbf{r}_2\] is a finite-codimensional submodule of \[\textstyle T(\ind_{KZ}^G \QUOTIENTI(\alpha)).\]\\[-22pt]
                \item\label{CORO1.6b} \m{\Big(\PADICVALUATION(a)-\alpha = v = \PADICVALUATION(\TEMPORARYD') \LEQSLANT \PADICVALUATION(C_{-1}) \annd \veeC\not\in \ZZ_p^\times \annd C_0 \in \ZZ_p^\times,\,\mathbf{r}_3\Big)} for
                \[\textstyle \mathbf{r}_3 = \left(T+\vveeC\begin{Lmatrix} 1 & 0 \\ 0 & p\end{Lmatrix}-\frac{C_0^{-1}\TEMPORARYD'}{ap^{-\alpha}}\right)(\ind_{KZ}^G \QUOTIENTI(\alpha)),\]
                where
                \[\textstyle \vveeC = \VERIFYIF{r\equiv_{p-1}2\alpha}\left((-1)^{\alpha}\sum_{l=0}^\alpha C_0^{-1}C_l \binom{r-\alpha+l}{s-\alpha}-1\right).\]\\[-22pt]
                                \item\label{CORO1.7b} \m{\Big(\PADICVALUATION(a)-\alpha = v = \PADICVALUATION(C_{-1}) < \PADICVALUATION(\TEMPORARYD') \annd  C_0 \not\in \ZZ_p^\times \annd \LESS{2\alpha-r}> 0,\,\mathbf{r}_4\Big)} for
                \[\textstyle \mathbf{r}_4 = \left(T+\frac{\veeC ap^{-\alpha}}{C_{-1}}\right)(\ind_{KZ}^G \SUBMODULE(\alpha)).\]\\[-22pt]
                                                \item\label{CORO1.8b} \m{\Big(\PADICVALUATION(a)-\alpha = v = \PADICVALUATION(C_{-1}) < \PADICVALUATION(\TEMPORARYD') \annd  \veeC \in \ZZ_p^\times,\,\ind_{KZ}^G \QUOTIENTI(\alpha)\Big)}\NEWLINEA                         \item\label{CORO1.9b} \m{\Big(\PADICVALUATION(a)-\alpha = v = \PADICVALUATION(C_{-1}) < \PADICVALUATION(\TEMPORARYD') \annd  C_0 \in \ZZ_p^\times,\, \mathbf{r}_5 \Big)},
                                                 where \[\textstyle\mathbf{r}_5\] is a finite-codimensional submodule of \[\textstyle T(\ind_{KZ}^G \QUOTIENTI(\alpha)).\]
\end{enumerate}
\end{corollary}

\begin{proof}{Proof. } %
(\ref{CORO1.1b}, \ref{CORO1.2b}, \ref{CORO1.3b}, \ref{CORO1.4b}, \ref{CORO1.5b}) The proofs of these parts are nearly identical to the proofs of the corresponding parts of corollary~\ref{CORO1}. 

(\ref{CORO1.6b}) 
The proof is similar to the proof of~(\ref{CORO1.5b}), the only difference being that the valuation of \m{\TEMPORARYD'} is the same as the valuation of the coefficient of \m{H}. To be more specific, we divide \m{\gamma} by  \m{C_0}, the term ``\m{T}'' comes from the expression for \m{H} given in lemma~\ref{LEMMATHM}, the term ``\m{\vveeC\begin{Lmatrix} 1 & 0 \\ 0 & p\end{Lmatrix}}'' comes from \[\textstyle \VERIFYIF{r\equiv_{p-1}2\alpha} (C_0^{-1}B-1) \begin{Lmatrix} 0 & 1 \\ p & 0\end{Lmatrix},\] the reason there is no term ``\m{A\begin{Lmatrix}
p & 0 \\ 0 & 1
\end{Lmatrix}}'' is because \m{A =\veeC = \BIGO{p}}, and the term ``\m{-\frac{C_0^{-1}\TEMPORARYD'}{ap^{-\alpha}}}'' comes from the first line of the formula for \m{\gamma} given in~\eqref{EQIMPX}.

(\ref{CORO1.7b}) 
As in the previous parts we can deduce that \m{\SCR I_a} contains \[\textstyle L \defeq \frac{\veeC ap^{-\alpha}}{C_{-1}} \begin{Lmatrix}
p & 0 \\0 & 1
\end{Lmatrix} \sqbrackets{}{KZ}{\QQp}{\theta^\alpha (y^{r'} - x^{p-1}y^{r'-p+1})} + 1 \sqbrackets{}{KZ}{\QQp}{\theta^\alpha y^{r'}} + L',\]
where \m{r' = r-\alpha(p+1)} and \m{L'} reduces modulo \m{\maxideal} to a trivial element of \m{\SUBMODULE(\alpha)}. The reduction  modulo \m{\maxideal} of the element
\m{\textstyle \sum_{\MU\in \FF_p} \begin{Smatrix} 1 & 0 \\ [\MU] & p \end{Smatrix}L} generates \m{\mathbf{r}_4}.

(\ref{CORO1.8b}, \ref{CORO1.9b}) The proofs of these parts are similar to the proofs of (\ref{CORO1.4b}, \ref{CORO1.5b}).
\end{proof}

\section{Proof of theorem~\ref{MAINTHM0a}}
\label{proof1}

The proof of theorem~\ref{MAINTHM0a} is based on the approach outlined in~\cite{BG09}, and  roughly consists of finding enough elements in \m{\SCR I_{a}}, consequently eliminating enough subquotients of \m{\ind_{KZ}^G\SIGMA_r}, and using that information to find \m{\overline\Theta_{k,a}}.

Throughout this section we assume that \[\textstyle r = s + \beta(p-1) + u_0p^t + \BIGO{p^{t+1}}\] for some \m{\beta \in \{0,\ldots,p-1\}}  and \m{u_0 \in \ZZ_p^\times} and \m{t \in \ZZ_{>0}}. Let us write \m{\TEMPORARYJ = u_0p^t}. Recall also that we assume 
\m{\TEMPORARYI-1<\PADICVALUATION(a)<\TEMPORARYI} for some  \m{\TEMPORARYI\in\{1,\ldots,\frac{p-1}{2}\}}, that \m{s\in\{2\TEMPORARYI,\ldots,p-2\}}, and that \m{k>p^{100}} (and consequently \m{r > p^{99}}).

Let us first show the equivalence between theorem~\ref{MAINTHM0a} and the union of the  following two propositions.

\begin{proposition}\label{PROP1.1}
If \m{\textstyle k\in \RRRRRR_{0}} then any infinite-dimensional factor of \m{\overline \Theta_{k,a}} is a quotient of \m{\FILT N_{\TEMPORARYI-1}}. If \m{\textstyle k\in \RRRRRR_{\beta,\jj}} then none of the infinite-dimensional factors of \m{\overline \Theta_{k,a}} are quotients of a representation in the set
\[\textstyle\left\{\FILT N_0,\ldots,\FILT N_{\max\{\TEMPORARYI-\jj-1,\beta\}-1}\right\}.\]
\end{proposition}

\begin{proposition}\label{PROP2.1}
If \m{\textstyle k\in \RRRRRR_{\beta,\jj}} then none of the infinite-dimensional factors of \m{\overline \Theta_{k,a}} are quotients of a representation in the set \[\textstyle\left\{\FILT N_{\max\{\TEMPORARYI-\jj-1,\beta\}+1},\ldots,\FILT N_{\TEMPORARYI-1}\right\}.\]
\end{proposition}

\begin{proof}{Proof that theorem~\ref{MAINTHM0a} is equivalent to propositions~\ref{PROP1.1} +~\ref{PROP2.1}. } 
First let us assume that propositions~\ref{PROP1.1} and~\ref{PROP2.1} are true. Together they imply that any infinite-dimensional factor of \m{\overline \Theta_{k,a}} is a quotient of \m{\FILT N_\alpha}, where \m{\alpha = \TEMPORARYI-1} if 
\m{\textstyle k\in \RRRRRR_{0}} and \m{\alpha = \max\{\TEMPORARYI-\jj-1,\beta\}} if \m{\textstyle k\in \RRRRRR_{\beta,\jj}}. 
The classifications given by theorem~2 in~\cite{Ars1} and theorem~2.7.1 in~\cite{Bre03a} imply that if \m{\overline \Theta_{k,a}} is reducible then it must have exactly two infinite-dimensional factors. There may be an additional one-dimensional factor, a twist of the Steinberg representation. Suppose that the infinite-dimensional factors are quotients of \m{\ind_{KZ}^G (\sigma_b(b'))} and \m{\ind_{KZ}^G (\sigma_d(d'))}, respectively. By theorem~2 in~\cite{Ars1} we must have
\[\textstyle b'-d'\equiv_{p-1} d+1 \text{ and } b+d\equiv_{p-1}-2.\]
In particular, \m{\ind_{KZ}^G (\sigma_b(b'))} and \m{\ind_{KZ}^G (\sigma_d(d'))} cannot be the representations \m{\ind_{KZ}^G\SUBMODULE(\alpha)} and \m{\ind_{KZ}^G\QUOTIENTI(\alpha)}, as that would imply that \[\textstyle   2\alpha-r\equiv_{p-1}2\alpha-r+1.\]
Similarly, the two representations cannot be two copies of
\m{\ind_{KZ}^G\SUBMODULE(\alpha)}, as that would imply that
\[\textstyle  2\alpha\equiv_{p-1}s+1,\]
a contradiction since \m{2\alpha\in \{2,\ldots,2\TEMPORARYI-2\}} and \m{s \in \{2\TEMPORARYI,\ldots,p-1\}}.  And, the two representations cannot be two copies of \m{\ind_{KZ}^G\QUOTIENTI(\alpha)}, as that would imply that 
 \[\textstyle   2\alpha\equiv_{p-1}s-1,\]
 which is similarly a contradiction. Thus we can conclude that
 \m{\overline \Theta_{k,a}} must be irreducible, and the classifications given by theorem~2 in~\cite{Ars1} and theorem~2.7.1 in~\cite{Bre03a} imply that the only possible quotient of \m{\FILT N_\alpha} that \m{\overline \Theta_{k,a}} can be is \m{\BIRREDUCIBLE(b_\alpha)}. This implies theorem~\ref{MAINTHM0a}.
 
Conversely, if theorem~\ref{MAINTHM0a} is true, then the fact that \m{\overline \Theta_{k,a} \cong\BIRREDUCIBLE(b_\alpha)} implies that \m{\overline \Theta_{k,a}} is irreducible and not a quotient of  a representation in the set \[\textstyle\left\{\FILT N_{0},\ldots,\FILT N_{\TEMPORARYI-1}\right\} \backslash \left\{N_\alpha\right\}.\]
\end{proof}

Proposition~\ref{PROP1.1} is essentially the main result of~\cite{Ars1}. Let us now prove proposition~\ref{PROP2.1}.

\emph{Proof of proposition~\ref{PROP2.1}. }
Let \m{\alpha \in \{0,\ldots,\TEMPORARYI-2\}} and let us consider \m{\FILT N_\alpha}. The task is to show that if \m{\alpha >\max\{\TEMPORARYI-\jj-1,\beta\}} then none of the infinite-dimensional factors of \m{\overline \Theta_{k,a}} are quotients of \m{\FILT N_\alpha}.
Note that the condition on \m{\alpha} implies both \m{\alpha > \beta} and \[\textstyle \jj \GEQSLANT \TEMPORARYI -\alpha > \PADICVALUATION(a)-\alpha.\] Let us apply part~(\ref{CORO1.3})  of corollary~\ref{CORO1} with \m{v} chosen arbitrarily in the open interval \m{(\PADICVALUATION(a)-\alpha,\jj)} and
\begin{align*} 
 C_j = \left\{\begin{array}{rl}  
  0 & \text{if } j  \in \{-1,0\}, \Cjnewline
 (-1)^{\alpha-j}\binom{s-\alpha+1}{\alpha-j} + pC_j^\ast & \text{if } j \in \{1,\ldots,\alpha\},
 \end{array}\right.\end{align*}
 for some constants \m{C_1^\ast,\ldots,C_\alpha^\ast} yet to be chosen. %
  We need to show that the constants \m{\{C_j\}} are suitable, \IE that the conditions of corollary~\ref{CORO1} are satisfied. Clearly  \[\textstyle \PADICVALUATION(a)-\alpha < v\LEQSLANT\PADICVALUATION(C_{-1})\] and \m{C_0 = \BIGO{p}}. Moreover,
 \begin{align*}\textstyle \sum_{l=1}^\alpha C_l \binom{r-\alpha+l}{l} & \textstyle = \sum_{l=1}^\alpha (-1)^{\alpha-l}\binom{s-\alpha+1}{\alpha-l}\binom{s-\alpha-\beta+l}{l} + \BIGO{p} \\ & \textstyle = (-1)^\alpha\binom{\beta}{\alpha} + (-1)^{\alpha+1} \binom{s-\alpha+1}{\alpha} + \BIGO{p}
 \\ & \textstyle\null = (-1)^{\alpha+1} \binom{s-\alpha+1}{\alpha} + \BIGO{p}.\end{align*}
 The third equality follows from the fact that \m{\alpha > \beta}. Since \[\textstyle p+\alpha-1>s>2\alpha-2\] for \m{\alpha>0}, and \m{\binom{s-\alpha+1}{\alpha}=1} for \m{\alpha=0}, it follows that   \[\textstyle \sum_{l=1}^\alpha C_l \binom{r-\alpha+l}{l} \in \ZZ_p^\times.\] Thus we only need to verify the most delicate condition, that \[\textstyle v\LEQSLANT \PADICVALUATION(\TEMPORARYD_{w}(D_\DARKCIRC))\] for \m{0\LEQSLANT w<2\TEMPORARYI-\alpha}.
By~\eqref{ceqn1} and~\eqref{ceqn6}, if
\begin{align*}
    & \textstyle L_1(r) \defeq \sum_{i>0} \binom{r-\alpha+j}{i(p-1)+j}\binom{i(p-1)}{w}
\end{align*} then \m{L_1(r) = L_1(s+\beta(p-1)) + \BIGO{\TEMPORARYJ}} for \m{0\LEQSLANT w<2\TEMPORARYI-\alpha}. So in order to verify the last condition it is enough to show that
    \begin{align}\label{eqn23}\textstyle L_2(s)\defeq\sum_{j=1}^\alpha \left((-1)^{\alpha-j} \binom{s-\alpha+1}{\alpha-j} + pC_j^\ast\right) \binom{s+\beta(p-1)-\alpha+j}{i(p-1)+j} = \BIGO{p^v}\end{align}
    for all \m{i\in\{1,\ldots,\beta\}}. %
     We have \[\textstyle \binom{s+\beta(p-1)-\alpha+j}{i(p-1)+j} = \binom{\beta}{i}\binom{s-\alpha-\beta+j}{j-i}+\BIGO{p}.\] We also have \begin{align*} & \textstyle \textstyle  %
\sum_{j=1}^\alpha (-1)^{\alpha-j} \binom{s-\alpha+1}{\alpha-j}\binom{s-\alpha-\beta+j}{j-i} \\ & \textstyle \null\qquad\tempheightaer = %
\sum_{j=1}^\alpha (-1)^{\alpha-i} \binom{s-\alpha+1}{\alpha-j}\binom{\alpha+\beta-s-i-1}{j-i}
    \\ & \textstyle \null\qquad\tempheightaer = (-1)^{\alpha-i} %
    \left(- \binom{s-\alpha+1}{\alpha}\binom{\alpha+\beta-s-i-1}{-i} + \sum_{j}  \binom{s-\alpha+1}{\alpha-j}\binom{\alpha+\beta-s-i-1}{j-i}\right)
        \\ & \textstyle \null\qquad\tempheightaer = (-1)^{\alpha-i} %
        \left(- \VERIFYIF{i=0}\binom{s-\alpha+1}{\alpha} +  \binom{\beta-i}{\alpha-i}\right) = 0.\end{align*}
The second equality follows from Vandermonde's convolution formula.  The third equality follows from the assumptions that \m{i>0} and \m{\alpha>\beta}. So~\eqref{eqn23} is true modulo~\m{p}, and we can transform~\eqref{eqn23} into the matrix equation
     \[\textstyle \left(\binom{s+\beta(p-1)-\alpha+j}{i(p-1)+j}\right)_{1\LEQSLANT i\LEQSLANT \beta, \, 1\LEQSLANT j \LEQSLANT \alpha} (C_1^\ast,\ldots,C_\alpha^\ast)^T = (w_1,\ldots,w_\beta)^T\]
     for some \m{w_1,\ldots,w_\beta \in \ZZ_p}. This matrix equation always has a solution since the left \m{\beta\times \beta} submatrix of the reduction modulo~\m{p} of the matrix
     \[\textstyle \left(\binom{s+\beta(p-1)-\alpha+j}{i(p-1)+j}\right)_{1\LEQSLANT i\LEQSLANT \beta, \, 1\LEQSLANT j \LEQSLANT \alpha}\]
     is upper triangular with units on the diagonal. Therefore we can indeed always choose the constants 
     \m{C_1^\ast,\ldots,C_\alpha^\ast} in a way that \m{v\LEQSLANT \PADICVALUATION(\TEMPORARYD_{w}(D_\DARKCIRC))} for \m{0\LEQSLANT w<2\TEMPORARYI-\alpha}. Then all conditions of part~(\ref{CORO1.3})  of corollary~\ref{CORO1} are satisfied, which concludes the proof of proposition~\ref{PROP2.1}.
\hspace{\stretch{1}}     {\BLACKSQUARE}

Propositions~\ref{PROP1.1} and~\ref{PROP2.1} are already sufficient to compute \m{\overline \Theta_{k,a}}, since they imply that \m{\overline \Theta_{k,a}} must be a quotient of \m{\FILT N_\alpha}, where \m{\alpha = \TEMPORARYI-1} if 
\m{\textstyle k\in \RRRRRR_{0}} and \m{\alpha = \max\{\TEMPORARYI-\jj-1,\beta\}} if \m{\textstyle k\in \RRRRRR_{\beta,\jj}}. Let us show that in fact  the surjective map \[\textstyle \FILT N_{\alpha} \xtwoheadrightarrow{\hspace{\SIZEEE} \hphantom{\cong}\hspace{\SIZEEE}} \overline\Theta_{k,a}\] factors through \m{\ind_{KZ}^G\QUOTIENTI(\alpha)}.

\begin{proposition}\label{PROP1.1bb}
If \m{k\in \RRRRRR_{0}} then the surjective map \[\textstyle \FILT N_{\TEMPORARYI-1} \xtwoheadrightarrow{\hspace{\SIZEEE} \hphantom{\cong}\hspace{\SIZEEE}} \overline\Theta_{k,a}\] factors through \m{\ind_{KZ}^G\QUOTIENTI(\TEMPORARYI-1)}.
\end{proposition}

\begin{proof}{Proof. }
Let \m{\alpha=\TEMPORARYI-1}. Note that since \m{k\in \RRRRRR_{0}} we have \[\textstyle\beta\in \{\TEMPORARYI-1,\ldots,p-1\}.\] Let us apply part~(\ref{CORO1.2}) of corollary~\ref{CORO1} with \m{v=0} and some constants \[\textstyle C_{-1},C_0,\ldots,C_\alpha\] such that \m{C_{-1}=\binom{\beta}{\alpha}} and \m{C_0=0}. The conditions that need to be satisfied in order for the corollary to be applicable are \m{\PADICVALUATION(\TEMPORARYD_{w}(D_\DARKCIRC))>0} for \m{0\LEQSLANT w<\alpha}, and \m{\PADICVALUATION(\TEMPORARYD')>0}. Let us consider the matrix \m{A = (A_{w,j})_{0\LEQSLANT w,j \LEQSLANT\alpha}} that has integer entries
\begin{align*}
     A_{w,j}  & \textstyle = \sum_{i>0} \binom{r-\alpha+j}{i(p-1)+j}\binom{i(p-1)}{w}.
\end{align*}
Then \m{\PADICVALUATION(\TEMPORARYD')>0} is equivalent to \m{\TEMPORARYD_{\alpha}(D_\DARKCIRC)=C_{-1}+\BIGO{p} = \binom{\beta}{\alpha}+\BIGO{p}}, so the two equations we want to show are equivalent to
\[\textstyle A (0,C_1,\ldots,C_\alpha)^T = \left(-\binom{\beta}{\alpha},0,\ldots,0,\binom{\beta}{\alpha}\right)^T + \BIGO{p}.\]%
By~\eqref{ceqn6} we have
\begin{align*}
     \overline A_{w,j}  & \textstyle = F_{w,j}(r,s),
\end{align*}
where \m{F_{w,j}(\polyr,\polys) \in \FF_p[\polyr,\polys]} is the polynomial defined in lemma~\ref{lemmaComb555}. Then the conclusion of that lemma is that
\[\textstyle \overline A (0,C_1,\ldots,C_\alpha)^T = \left(-\binom{\beta}{\alpha},0,\ldots,0,\binom{\beta}{\alpha}\right)^T \]
with
\[\textstyle C_j = (-1)^{j} \binom{s-\alpha+1}{\alpha-j}\] for \m{j \in \{1,\ldots,\alpha\}}.
Thus these choices for \m{C_{-1},C_0,\ldots,C_\alpha} are suitable, and we can apply part~(\ref{CORO1.2})  of corollary~\ref{CORO1} with \m{v = 0} and conclude that  \m{\ind_{KZ}^G\SUBMODULE(\alpha)} is trivial modulo \m{\SCR I_a}.
\end{proof}

\begin{proposition}\label{PROP1.1b}
If \m{k\in \RRRRRR_{\beta,\jj}} then the surjective map \[\textstyle \FILT N_{\max\{\TEMPORARYI-\jj-1,\beta\}} \xtwoheadrightarrow{\hspace{\SIZEEE} \hphantom{\cong}\hspace{\SIZEEE}} \overline\Theta_{k,a}\] factors through \m{\ind_{KZ}^G\QUOTIENTI(\max\{\TEMPORARYI-\jj-1,\beta\})}.
\end{proposition}

\begin{proof}{Proof. }
Let \m{\alpha=\max\{\TEMPORARYI-\jj-1,\beta\}}.
Let us apply part~(\ref{CORO1.2})  of corollary~\ref{CORO1} with \m{v=\jj} and
 \begin{align*} 
 C_j = \left\{\begin{array}{rl}  \TEMPORARYJ & \text{if } j = -1,\Cjnewline
 0 & \text{if } j = 0,\Cjnewline
 (-1)^{\alpha+\beta+j} \alpha \binom{\alpha-1}{\beta} \binom{s-\alpha+1}{\alpha-j} & \text{if } j \in \{1,\ldots,\alpha\}.
 \end{array}\right.\end{align*}
 We need to show that the constants \m{\{C_j\}} are suitable, \IE that the conditions of corollary~\ref{CORO1} are satisfied.
 Clearly \m{\PADICVALUATION(C_{-1})=\jj < \PADICVALUATION(a)-\alpha}. We also need to show that \m{\jj<\PADICVALUATION(\TEMPORARYD_{w}(D_\DARKCIRC))} for \m{0\LEQSLANT w<\alpha} and \m{\jj < \PADICVALUATION(\TEMPORARYD')} and \m{\jj\LEQSLANT\PADICVALUATION(\TEMPORARYD_{w}(D_\DARKCIRC))} for \m{\alpha\LEQSLANT w<2\TEMPORARYI-\alpha}. Let us consider the matrix \m{A = (A_{w,j})_{0\LEQSLANT w,j \LEQSLANT\alpha}} that has integer entries
\begin{align*}
     A_{w,j}  & \textstyle = \sum_{i>0} \binom{r-\alpha+j}{i(p-1)+j}\binom{i(p-1)}{w}.
\end{align*}
 If we consider the approximation claim in the proof of proposition~\ref{PROP1.1} and multiply the first column by \m{p}, we get
\[\textstyle A = S + \TEMPORARYJ N + \BIGO{\TEMPORARYJ p},\] where
\begin{align*}
    S_{w,j} & = \textstyle \sum_{i=1}^\beta \binom{s+\beta(p-1)-\alpha+j}{i(p-1)+j}\binom{i(p-1)}{w},\\
    N_{w,j} & = \textstyle \sum_v (-1)^{w-v} \binom{j+w-v-1}{w-v} \binom{s+\beta(p-1)-\alpha+j}{v}^\partial \sum_{i=0}^\beta \binom{s+\beta(p-1)-\alpha+j-v}{i(p-1)+j-v}\\ & \textstyle\null\qquad -\VERIFYIF{w=0} \binom{s+\beta(p-1)-\alpha+j}{j}^\partial.
\end{align*}
Exactly as in the proof of proposition~\ref{PROP1.1} we can deduce the three conditions we need to show as long as \begin{align*} \textstyle S(C_0,\ldots,C_\alpha)^T & \textstyle  = 0, \\ 
\textstyle N (C_0,\ldots,C_\alpha)^T & \textstyle = (-C_{-1}\TEMPORARYJ^{-1},0,\ldots,0,C_{-1}\TEMPORARYJ^{-1})^T + Sv + \BIGO{p} \text{ for some }v.\end{align*}%
 Let
\m{B = B_\alpha} be the 
\m{(\alpha+1)\times(\alpha+1)} matrix defined in lemma~\ref{lemmaComb2}. That lemma implies that 
\m{B} encodes precisely the row operations that transform \m{S} into a matrix with zeros outside the rows indexed \m{1,\ldots,\beta} and such that \[\textstyle (BS)_{w,j} = p^{-\VERIFYIF{j=0}}\binom{s+\beta(p-1)-\alpha+j}{w(p-1)+j}\]
when \m{w \in\{1,\ldots,\beta\}}. Moreover,
\[\textstyle {B}_{i,w} = \VERIFYIF{(i,w)=(0,0)} + \sum_{l=1}^\alpha (-1)^{i} \binom{l}{i} \binom{l-1}{l-w} + \BIGO{p}.\]
By using this formula we can compute that
 \[\textstyle B (-1,0,\ldots,0,1)^T = \left(0,-\binom{\alpha}{1},\ldots,(-1)^\alpha \binom{\alpha}{\alpha}\right)^T + \BIGO{p},\]
 and therefore if \m{\overline{R}} is the \m{\alpha\times\alpha} matrix over \m{\FF_p} obtained from \m{\overline{BN}} by replacing the rows indexed \m{1,\ldots,\beta} with the corresponding rows of \m{\overline{BS}} and then discarding the zeroth row and the zeroth column, the condition that needs to be satisfied is equivalent to the claim that
 \begin{align}\label{eq34566}\textstyle 
 \overline R (C_1,\ldots,C_\alpha)^T = \left(-(1-\VERIFYIF{1\LEQSLANT\beta})\binom{\alpha}{1},\ldots,(-1)^\alpha (1-\VERIFYIF{\alpha\LEQSLANT\beta}) \binom{\alpha}{\alpha}\right)^T.\end{align}
 The matrix \m{\overline R} is the lower right  \m{\alpha\times\alpha} submatrix of the matrix \m{\overline Q} defined in the proof of proposition~\ref{PROP1.1}, where we compute that
 \begin{align*} \textstyle \overline R_{i-1,j-1} & \textstyle =\binom{s-\alpha-\beta+j}{j-i} \cdot \left\{\begin{array}{rl}
    \binom{\beta}{i}  & \text{if } i \in \{1,\ldots,\beta\}, \Cjnewline
    -\binom{\beta}{i}^\partial  & \text{otherwise,}
 \end{array}\right.\end{align*}
 for \m{i,j \in \{1,\ldots,\alpha\}}. We have
 \begin{align*}
     & 
     \textstyle \sum_{j=1}^\alpha (-1)^{\alpha+\beta+j} \alpha \binom{\alpha-1}{\beta} \binom{s-\alpha+1}{\alpha-j} \binom{s-\alpha-\beta+j}{j-i} \\& \textstyle \null\qquad =\sum_{j} (-1)^{\alpha+\beta+i} \alpha \binom{\alpha-1}{\beta} \binom{s-\alpha+1}{\alpha-j} \binom{\alpha+\beta-s-i-1}{j-i}
     \\& \textstyle \null\qquad=(-1)^{\alpha+\beta+i} \alpha \binom{\alpha-1}{\beta} \binom{\beta-i}{\alpha-i}.
\end{align*}
If \m{i \in \{1,\ldots,\beta\}} then the last expression is zero, and if \m{i\in \{\beta+1,\ldots,\alpha\}}
 then it is
 \[\textstyle (-1)^{\beta} \alpha \binom{\alpha-1}{\beta} \binom{\alpha-\beta-1}{i-\beta-1}%
 = (-1)^\beta i\binom{i-1}{\beta} \binom{\alpha}{i}   = (-1)^{i} \binom{\alpha}{i} \cdot \left(-\binom{\beta}{i}^\partial\right)^{-1},\]
 which implies~\eqref{eq34566}.
Consequently we can apply part~(\ref{CORO1.2})  of corollary~\ref{CORO1} with \m{v = \jj} and conclude that  \m{\ind_{KZ}^G\SUBMODULE(\alpha)} is trivial modulo \m{\SCR I_a}.
\end{proof}

\section{Proof of theorem~\ref{MAINTHM00a}}

The proof of theorem~\ref{MAINTHM00a} is very similar to the proof of theorem~\ref{MAINTHM0a}. The major difference is that we apply corollary~\ref{CORO2} instead of corollary~\ref{CORO1} since \m{\PADICVALUATION(a)} is an integer, which means that \m{\overline{\Theta}_{k,a}} is reducible in some cases.

We make the same assumptions as in section~\ref{proof1}.
We assume that \[\textstyle r = s + \beta(p-1) + u_0p^t + \BIGO{p^{t+1}}\] for some \m{\beta \in \{0,\ldots,p-1\}}  and \m{u_0 \in \ZZ_p^\times} and \m{t \in \ZZ_{>0}}, that \m{\TEMPORARYJ = u_0p^t}, that
\m{\PADICVALUATION(a)=\TEMPORARYI-1} for some  \m{\TEMPORARYI\in\{1,\ldots,\frac{p-1}{2}\}}, that \m{s\in\{2\TEMPORARYI,\ldots,p-2\}}, and that \m{k>p^{100}} (and consequently \m{r > p^{99}}).

Let \m{\mu = \lambda} if \m{k \in \RRRRRR_0} and \m{\mu = \lambda_{\beta,t}} if \m{k\in \RRRRRR_{\beta,t}} (in the notation of theorem~\ref{MAINTHM00}). Let  \m{\llll=\TEMPORARYI-1} if \m{k \in \RRRRRR_0} and \m{\llll = \max\{\TEMPORARYI-\jj-1,\beta\}} if \m{k\in \RRRRRR_{\beta,\jj}}. Let us first show the equivalence between theorem~\ref{MAINTHM00a} and the following proposition.

\begin{proposition}\label{PROP3.1}
Let either \m{k \in \RRRRRR_0}, or \m{k \in \RRRRRR_{\beta,t}} and \m{t<\TEMPORARYI-\beta-1}.
\begin{enumerate} 
\item\label{itemprop1.2part1}\m{ (T-\mu^{-1}) (\ind_{KZ}^G\QUOTIENTI(\llll-1))} is trivial modulo~\m{\SCR I_a}. %
\\
\item\label{itemprop1.2part2}\m{(T-\mu) (\ind_{KZ}^G\SUBMODULE(\llll))} is trivial modulo~\m{\SCR I_a}. %
\\
\item\label{itemprop1.2part0}\m{\ind_{KZ}^G\SUBMODULE(\llll-1)} is trivial modulo~\m{\SCR I_a}. %
\\
\item\label{itemprop1.2part3} \m{\ind_{KZ}^G\QUOTIENTI(\llll)} is trivial modulo~\m{\SCR I_a}.
\end{enumerate}
\end{proposition}

\begin{proof}{Proof that theorem~\ref{MAINTHM00a} is equivalent to proposition~\ref{PROP3.1}. }
First let us assume that proposition~\ref{PROP3.1} is true.
In the setting of theorem~\ref{MAINTHM0a}, propositions~\ref{PROP1.1}, \ref{PROP2.1}, \ref{PROP1.1bb}, and~\ref{PROP1.1b} show  that if \m{k\in \RRRRRR_{0}} then  \m{\overline \Theta_{k,a}} is a quotient of \m{\ind_{KZ}^G\QUOTIENTI(\TEMPORARYI-1)}, and if \m{k\in \RRRRRR_{\beta,\jj}} then \m{\overline \Theta_{k,a}} is a quotient of \m{\ind_{KZ}^G\QUOTIENTI(\max\{\TEMPORARYI-\jj-1,\beta\})}. Their proofs are based on corollary~\ref{CORO1}. They amount to considering the element of \m{\Img(T-a)} coming from equation~\ref{EQIMPX} in lemma~\ref{LEMMATHM}, and noting that the term with dominant valuation is either \m{H} or
\[\textstyle (\TEMPORARYD'+C_{-1}) \sqbrackets{}{KZ}{\QQp}{\theta^\alpha x^{p-1}y^{r-\alpha(p+1)-p+1}}  + C_{-1}  \sqbrackets{}{KZ}{\QQp}{\theta^n x^{\alpha-n}y^{r-np-\alpha}},\] depending on how \m{\jj} compares to \m{\PADICVALUATION(a)-\alpha}.
In the setting of theorem~\ref{MAINTHM00a} we can apply the analogous corollary~\ref{CORO2} to conclude that any infinite-dimensional factor of \m{\overline \Theta_{k,a}} must be a quotient of one of
\[\textstyle \ind_{KZ}^G\SUBMODULE(\llll-1),\, \ind_{KZ}^G\QUOTIENTI(\llll-1),\, \ind_{KZ}^G\SUBMODULE(\llll),\, \ind_{KZ}^G\QUOTIENTI(\llll),\]
where for convenience we define \m{\SUBMODULE(-1)} and \m{\QUOTIENTI(-1)} to be the trivial representation. The key reason why the proofs of propositions~\ref{PROP1.1} and~\ref{PROP2.1} copy verbatim to prove this is that outside of these subquotients the valuations \m{\jj} and \m{\PADICVALUATION(a)-\alpha} never match, so again exactly one of the two aforementioned terms is dominant. The only subtlety when copying the proofs of propositions~\ref{PROP1.1} and~\ref{PROP2.1} is that we do not know whether \m{\SCR I_a} contains \m{1\sqbrackets{}{KZ}{\QQp}{x^{\TEMPORARYI-1}y^{r-\TEMPORARYI+1}}}. This ultimately does not present a problem since when working with \m{\FILT N_{\TEMPORARYI-1}} we always assume that \m{C_0 = 0}. As the proofs of propositions~\ref{PROP1.1} and~\ref{PROP2.1} work here nearly without modification except for replacing corollary~\ref{CORO1} with~\ref{CORO2}, we omit the full details of the arguments. Proposition~\ref{PROP3.1} then implies that any infinite-dimensional factor of  \m{\overline \Theta_{k,a}} must be a quotient of one of
\[\textstyle  \ind_{KZ}^G\QUOTIENTI(\llll-1)/(T-\mu^{-1}),\, \ind_{KZ}^G\SUBMODULE(\llll)/(T-\mu),\]
and together with theorem~2 in~\cite{Ars1} they completely determine~\m{\overline \Theta_{k,a}}.

The converse, that theorem~\ref{MAINTHM0a} implies proposition~\ref{PROP3.1}, is clear since  theorem~\ref{MAINTHM0a} completely determines \m{\overline \Theta_{k,a}}.
\end{proof}

\emph{Proof of proposition~\ref{PROP3.1}. }\label{subsectionMAINTHM00proof}

(\ref{itemprop1.2part1})
Suppose first that \m{k \in \RRRRRR_0}.
Let \m{c=(c_1,\ldots,c_\alpha)} be as in lemma~\ref{part4X5}, and let us make the substitutions \m{X=r-\alpha} and \m{Y=s-\alpha}.
We apply part~(\ref{CORO1.6b}) of corollary~\ref{CORO2} with \m{v=1}. We choose \m{C_{-1}=0} and \m{C_0 = 1} and \m{C_j = \frac{c_jp}{(s-\alpha)_\alpha}} for \m{j \in \{1,\ldots,\alpha\}}.
In the proof of proposition~\ref{PROP1.1} we show that for these constants we have
\begin{align*}
    \textstyle \TEMPORARYD' & \textstyle = \frac{(s-r)_{\TEMPORARYI-1}}{(s-\TEMPORARYI+2)_{\TEMPORARYI-1}}p + \BIGO{p^2}, \\
        1 & \textstyle\LEQSLANT \PADICVALUATION(\TEMPORARYD_{w}(D_\DARKCIRC))\text{ for } \alpha\LEQSLANT w<2\TEMPORARYI-\alpha, \\
    \textstyle 1 & \textstyle <\PADICVALUATION(\TEMPORARYD_{w}(D_\DARKCIRC)) \text{ for }0\LEQSLANT w<\alpha.
\end{align*}
Moreover, since \m{C_1,\ldots,C_\alpha=\BIGO{p}} we have \m{\veeC=\BIGO{p}}, and \[\textstyle\mu = \frac{(s-\TEMPORARYI+2)_{\TEMPORARYI-1}a}{(s-r)_{\TEMPORARYI-1} p^{\TEMPORARYI}} = \frac{a}{C_0^{-1}\TEMPORARYD' p^{\TEMPORARYI-2}}.\] Therefore the conditions needed to apply part~(\ref{CORO1.6b}) of corollary~\ref{CORO2} (\IE the three conditions \begin{align*}\textstyle  \textstyle v & \textstyle \LEQSLANT \PADICVALUATION(\TEMPORARYD_{\alpha}(D_\DARKCIRC)),
\\
v' \defeq \min\{\PADICVALUATION(a)-\alpha,v\}  & \textstyle \LEQSLANT \PADICVALUATION(\TEMPORARYD_{\TEMPORARYAAAA}(D_\DARKCIRC)) \text{ for } \alpha<\TEMPORARYAAAA<2\TEMPORARYI-\alpha,
\\ v'  & \textstyle <\PADICVALUATION(\TEMPORARYD_{\TEMPORARYAAAA}(D_\DARKCIRC)) \text{ for } 0\LEQSLANT w<\alpha,
\end{align*} in addition to the two extra conditions on \m{\veeC} and \m{C_0} in part~(\ref{CORO1.6b}) of the corollary) are satisfied and we can conclude that \[\textstyle (T-\mu^{-1}) (\ind_{KZ}^G\QUOTIENTI(\llll-1))\] is trivial modulo \m{\SCR I_a}. If \m{k \in \RRRRRR_{\beta,t}} and \m{t<\TEMPORARYI-\beta-1} then the argument is similar: we choose \m{v=t+1} and the same constants \m{\{C_j\}} as in the third bullet point (if \m{\beta=0}) or the fourth bullet point (if \m{\beta\in\{1,\ldots,\gamma-1\}}) of the proof of proposition~\ref{PROP1.1}. In the former case \[\textstyle \TEMPORARYD' = \frac{(s-r)_{\alpha+1}}{(s-\alpha)_{\alpha+1}}p + \BIGO{p^{t+2}},\] and in the latter case \[\textstyle \TEMPORARYD' = \epsilon Q_{0,0} p + \BIGO{p^{t+2}}.\] In both cases \[\textstyle \mu = \frac{a}{C_0^{-1}\TEMPORARYD' p^{\TEMPORARYI-t-2}}\] and the conditions needed to apply part~(\ref{CORO1.6b}) of corollary~\ref{CORO2} are satisfied, so again we can conclude that  \[\textstyle (T-\mu^{-1}) (\ind_{KZ}^G\QUOTIENTI(\llll-1))\] is trivial modulo \m{\SCR I_a}. %

(\ref{itemprop1.2part2})
Suppose first that \m{k \in \RRRRRR_0}. We use the constants 
\begin{align*} 
 C_j = \left\{\begin{array}{rl}  (-1)^{\TEMPORARYI-1}\binom{s-r}{\TEMPORARYI-1} & \text{if } j = -1,\Cjnewline
 0 & \text{if } j = 0,\Cjnewline
 (-1)^{\TEMPORARYI-j-1}\binom{s-\TEMPORARYI+2}{\TEMPORARYI-j-1} & \text{if } j \in \{1,\ldots,\alpha\}.
 \end{array}\right.\end{align*}
We can show just as in the proof of proposition~\ref{PROP1.1bb} that these constants satisfy all of the conditions needed to apply part~(\ref{CORO1.7b}) of corollary~\ref{CORO2} with \m{v=0}.  Moreover, we have
\begin{align*}\textstyle \veeC & \textstyle = \sum_{j=1}^{\TEMPORARYI-1} (-1)^{\TEMPORARYI-j-1}\binom{s-\TEMPORARYI+2}{\TEMPORARYI-j-1} \binom{r-\TEMPORARYI+j+1}{j} - (-1)^{\TEMPORARYI-1}\binom{s-r}{\TEMPORARYI-1}+ \BIGO{p}
\\ & \textstyle = (-1)^{\TEMPORARYI-1} \sum_{j=1}^{\TEMPORARYI-1} \binom{s-\TEMPORARYI+2}{\TEMPORARYI-j-1} \binom{\TEMPORARYI-r-2}{j} - (-1)^{\TEMPORARYI-1}\binom{s-r}{\TEMPORARYI-1}+ \BIGO{p}
\\ & \textstyle = (-1)^{\TEMPORARYI} \binom{s-\TEMPORARYI+2}{\TEMPORARYI-1} + \BIGO{p}.\end{align*}
Thus \[\textstyle-\frac{\veeC a p^{1-\TEMPORARYI}}{C_{-1}} = \frac{ \binom{s-\TEMPORARYI+2}{\TEMPORARYI-1}a}{\binom{s-r}{\TEMPORARYI-1}p^{\TEMPORARYI-1}} = \mu,\]
so  part~(\ref{CORO1.7b}) of corollary~\ref{CORO2} implies that  \[\textstyle (T-\mu) (\ind_{KZ}^G\SUBMODULE(\llll))\] is trivial modulo \m{\SCR I_a}. If \m{k \in \RRRRRR_{\beta,t}} and \m{t<\TEMPORARYI-\beta-1} then the argument is similar: we choose \m{v=t} and the constants 
 \begin{align*} 
 C_j = \left\{\begin{array}{rl}  \TEMPORARYJ & \text{if } j = -1,\Cjnewline
 0 & \text{if } j = 0,\Cjnewline
 (-1)^{\gamma+\beta+j} \gamma \binom{\gamma-1}{\beta} \binom{s-\gamma+1}{\gamma-j} & \text{if } j \in \{1,\ldots,\gamma\},
 \end{array}\right.\end{align*}
as in the proof of proposition~\ref{PROP1.1b}. Again all of the  conditions needed to apply part~(\ref{CORO1.7b}) of corollary~\ref{CORO2} with \m{v=t} are satisfied and
\begin{align*}\textstyle \veeC & \textstyle = (-1)^{\beta+\gamma} \gamma \binom{\gamma-1}{\beta} \sum_{j=1}^{\gamma} (-1)^j \binom{s-\gamma+1}{\gamma-j} \binom{r-\gamma+j}{j}  + \BIGO{p}
\\ & \textstyle = (-1)^{\beta+\gamma} \gamma \binom{\gamma-1}{\beta} \sum_{j=1}^{\gamma} \binom{s-\gamma+1}{\gamma-j} \binom{\gamma-r-1}{j}  + \BIGO{p}
\\ & \textstyle = (-1)^{\beta+\gamma} \gamma \binom{\gamma-1}{\beta} \left( \binom{s-r}{\gamma} - \binom{s-\gamma+1}{\gamma} \right)  + \BIGO{p}
\\ & \textstyle = (-1)^{\beta+\gamma+1} \gamma \binom{\gamma-1}{\beta} \binom{s-\gamma+1}{\gamma}   + \BIGO{p}.\end{align*}
The last equality follows from the fact that \m{\binom{s-r}{\gamma} = \BIGO{p}}. Thus \[\textstyle-\frac{\veeC a}{C_{-1}p^{\TEMPORARYI-\jj-1}} =  \frac{  (-1)^{\beta+\gamma} \gamma \binom{\gamma-1}{\beta} \binom{s-\gamma+1}{\gamma} a}{\epsilon p^{\TEMPORARYI-\jj-1}} = \mu,\]
so part~(\ref{CORO1.7b}) of corollary~\ref{CORO2} implies that   \[\textstyle (T-\mu) (\ind_{KZ}^G\SUBMODULE(\llll))\] is trivial modulo \m{\SCR I_a}.

(\ref{itemprop1.2part0})
This is very similar to part~(\ref{itemprop1.2part2}) of this proposition: if \m{k \in \RRRRRR_0} then we use part~(\ref{CORO1.2b})  of corollary~\ref{CORO2} just as in the proof of proposition~\ref{PROP1.1bb}, and if \m{k \in \RRRRRR_{\beta,t}} and \m{t<\TEMPORARYI-\beta-1} then we use part~(\ref{CORO1.2b})  of corollary~\ref{CORO2} just as in the proof of proposition~\ref{PROP1.1b}. We omit the full details.

(\ref{itemprop1.2part3}) We apply part~(\ref{CORO1.8b}) of corollary~\ref{CORO2} with the same constants as in the proof of part~(\ref{itemprop1.2part2}) of this proposition---since \m{\veeC \in \ZZ_p^\times} all of the necessary conditions are satisfied and we can conclude that  \m{\ind_{KZ}^G\QUOTIENTI(\llll)} is trivial modulo \m{\SCR I_a}.
\hspace{\stretch{1}}     {\BLACKSQUARE}

\end{document}